# LOCALLY ADAPTIVE ESTIMATION OF EVOLUTIONARY WAVELET SPECTRA[1]

By Sébastien Van Bellegem and Rainer von Sachs

*Université catholique de Louvain*

We introduce a wavelet-based model of local stationarity. This model enlarges the class of *locally stationary wavelet processes* and contains processes whose spectral density function may change very suddenly in time. A notion of *time-varying wavelet spectrum* is uniquely defined as a wavelet-type transform of the autocovariance function with respect to so-called *autocorrelation wavelets*. This leads to a natural representation of the autocovariance which is localized on scales. We propose a pointwise adaptive estimator of the time-varying spectrum. The behavior of the estimator studied in homogeneous and inhomogeneous regions of the wavelet spectrum.

**1. Introduction.** The spectral analysis of time series is a large field of great interest from both theoretical and practical viewpoints. The fundamental starting point for this analysis is the *Cramér representation*, stating that all zero-mean second order stationary processes $X_t$, $t \in \mathbb{Z}$, may be written

$$(1.1) \qquad X_t = \int_{[-\pi,\pi)} A(\omega) \exp(i\omega t) \, dZ(\omega), \qquad t \in \mathbb{Z},$$

where $A(\omega)$ is the amplitude of the process $X_t$ and $dZ(\omega)$ is an orthonormal increment process, that is, $\mathrm{E}(dZ(\omega)\overline{dZ(\mu)}) = d\omega \delta_0(\omega - \mu)$; see Brillinger (1975). Correspondingly, under mild conditions, the autocovariance function can be expressed as

$$c_X(\tau) = \int_{-\pi}^{\pi} f_X(\omega) \exp(i\omega\tau) \, d\omega,$$

Received June 2007; revised June 2007.

[1]Supported IAP research network nr P6/03 of the Belgian Government (Belgian Science Policy).

where $f_X$ is the *spectral density* of $X_t$.

There is not a unique way to relax the assumption of stationarity, that is, to define a second order process with a time-dependent spectrum. However, this modeling is a theoretical challenge which may be helpful in practice since many studies have shown that models with evolutionary spectra or time-varying parameters are necessary to explain some observed data, even over short periods of time. Examples may be found in numerous fields, such as economics [Swanson and White (1997), Los (2000)], biostatistics [Ombao et al. (2002)] and meteorology [Nason and Sapatinas (2002)], to name but a few.

Among the different possibilities for modeling nonstationary second order processes, we emphasize the approaches consisting of a modification of the Cramér representation (1.1). Different modifications of (1.1) are possible. First, we can replace the process $dZ(\omega)$ by a nonorthonormal process, leading to, for instance, the *harmonizable* processes [Lii and Rosenblatt (2002)]. A second possibility is to replace the amplitude function $A(\omega)$ by a time-varying version $A_t(\omega)$ and to assume a slow change of $A_t(\omega)$ over time. Such an approach is followed to define *oscillatory* processes [Priestley (1965)]. However, a major statistical drawback of the oscillatory processes is the intrinsic impossibility of constructing an asymptotic theory for consistency and inference. To overcome this problem, Dahlhaus (1997) introduced the class of *locally stationary processes*, in which the transfer function is rescaled in time. In this approach, a doubly-indexed process is defined as

$$(1.2) \quad X_{t,T} = \int_{[-\pi,\pi)} A\left(\frac{t}{T}, \omega\right) \exp(i\omega t) \, dZ(\omega), \qquad t = 0, \ldots, T-1, T > 0,$$

where the transfer function $A(z, \omega)$ is defined on $(0, 1) \times [-\pi, \pi)$. Dahlhaus (1997, 2000) investigated statistical inference for such processes, with a discussion on maximum likelihood, Whittle and least squares estimates, and showed that asymptotic results when $T$ tends to infinity can be considered. However, in this setting, letting $T$ tend to infinity does not have the usual meaning of "looking into the future," but means that we have, in the sample $X_{0,T}, \ldots, X_{T-1,T}$, more information about the local structure of $A(z, \omega)$. This formalism is analogous to nonparametric regression, for which "asymptotic" means an ideal knowledge about the local structure of the underlying curve.

In this article, we focus on a class of doubly-indexed locally stationary processes defined by replacing the harmonic system $\{\exp(i\omega t)\}$ in (1.2) by a wavelet system. In this way, we move from a *time-frequency* representation to a *time-scale* representation of the nonstationary process. Because wavelets systems are well localized in time and frequency, they appear more natural for modeling the time-varying spectra of nonstationary processes. As wavelets decompose the frequency domain into discrete scales, they offer a well-adapted system to achieve the trade-off resolution between time and frequency [Vidakovic (1999)].

The class of locally stationary wavelet processes studied in this article was initially introduced by Nason, von Sachs and Kroisandt (2000). Their definition of wavelet processes involves a time-varying amplitude which is smoothly varying and continuous as a function of time. An initial goal of this article is to extend this definition to the case of time-varying amplitudes with possibly *discontinuous* behavior in time. This introduces some technical difficulties to the proof of our results, but we believe the gain due to this extension to be crucial. Our new definition now includes more important examples of nonstationary processes. For instance, this extension of the definition is needed if we wish to model a nonstationary process built as a concatenation of different processes, such as the Haar processes defined in Nason et al. (2000). Moreover, wavelet processes can now be used for the analysis of intermittent phenomena, such as transients followed by regions of smooth behavior.

Our definition of wavelet processes is presented in Section 2, where we also define their *evolutionary spectrum*. This spectrum is a function of time and scales, and measures the power of the process at a particular time and scale. The main goal of the present article is to provide a pointwise adaptive estimation of the evolutionary spectrum. The estimation procedure follows the local adaptive method of Lepski (1990). The main difference with the latter is that we are now estimating a spectral density function, that is, the second order structure of correlated observations. Moreover, our statistical model is allowed to be nonstationary and the behavior of its evolutionary spectrum may be very inhomogeneous in time.

In Section 3, we present a preliminary estimator of the evolutionary spectrum and derive some useful properties that are needed in order to derive the adaptive estimator in Section 4. The behavior of this estimator is discussed for the two cases where the evolutionary wavelet spectrum is either regular or irregular near the point of estimation. These results explain the good performance of the algorithm in practice. Section 5 concludes with the result of a brief simulation study. All details and specific questions related to the practical implementation of our procedure have been considered in a separate paper [Van Bellegem and von Sachs (2004)], where a more exhaustive study of simulations and a real data analysis are provided.

Proofs and technical derivations are deferred to the appendices. Our estimator takes the form of a quadratic form of the increments, which are assumed to be Gaussian. Our estimator is the sum of a quadratic form of the increments that are assumed to be Gaussian and an additive, independent linear form of Gaussian variables. Thus, the main technical goal is to study the behavior of the (quadratic + linear) form of Gaussian variables. There exists a large body of results on quadratic forms of Gaussian variables. Recent developments include Rudzkis (1978), Neumann (1996), Laurent and Massart (2000), Spokoiny (2001), Comte (2001) and Dahlhaus and Polonik (2002). The exponential inequality proved in the latter reference is the starting point for some important results in the present article. On the other hand, in the appendices, we also present some original results on quadratic forms that are needed to prove our results.

**2. Locally stationary wavelet processes.** The wavelet system used to build locally stationary processes is a nondecimated system of compactly supported and discrete wavelets. We first briefly recall some points about this system of wavelets and then give a definition of the wavelet processes and wavelet spectra.

2.1. *Discrete nondecimated wavelet system.* The local functions used in the representation of LSW processes are a set of *discrete nondecimated wavelets* $\{\psi_{jk}, j = -1, -2, \ldots; k \in \mathbb{Z}\}$. We refer to Vidakovic (1999) for a review of wavelet theory and its applications in statistics, and to Nason and Silverman (1995) for a detailed introduction to the nondecimated wavelet transform. Let us simply recall that, in contrast to the discrete wavelet transform, the discrete nondecimated wavelets at all scales $j < 0$ can be shifted to any location defined by the finest resolution scale, determined by the observed data. As a consequence, this construction leads to an overcomplete system of the space of square summable sequences, $\ell^2(\mathbb{Z})$. The wavelets considered in this article are assumed to be compactly supported in time and we will denote by $\mathcal{L}_j$ the length of the support of $\psi_{j0}$, that is, $\mathcal{L}_j := |\operatorname{supp} \psi_{j0}|$. This automatically implies $|\operatorname{supp} \psi_{jk}| = \mathcal{L}_j = (2^{-j} - 1)(\mathcal{L}_{-1} - 1) + 1$ for all $j < 0$. Also, observe that, as in Nason et al. (2000), we departed from the usual wavelet numbering scheme. The data inhabit scale zero, and scale $-1$ is the scale which contains the finest resolution wavelet detail. Then, the support of the wavelet on the finest scale remains constant with respect to $T$.

For ease of presentation, recall the simplest discrete nondecimated system, called the *Haar system*, given by

$$\psi_{jk} = 2^{j/2} \mathbb{I}_{\{0,1,\ldots,2^{-j-1}-1\}}(k) - 2^{j/2} \mathbb{I}_{\{2^{-j-1},\ldots,2^{-j}-1\}}(k)$$
$$\text{for } j = -1, -2, \ldots \text{ and } k \in \mathbb{Z},$$

where $\mathbb{I}_{\mathcal{A}}(t)$ is 1 if $t \in \mathcal{A}$ and 0 otherwise. The shifted version of $\psi_{jk}$ is given by $\psi_{jk}(t) = \psi_{j,k-t}$ for all $k \in \mathbb{Z}$.

2.2. *The process and its evolutionary wavelet spectrum.* As we will note below, our definition of locally stationary wavelet processes differs from the original definition of Nason et al. (2000) as we only impose a total variation condition on the amplitudes instead of a Lipschitz condition. See also Fryźlewicz and Nason (2006) for a discussion of that definition.

DEFINITION 1. A sequence of doubly-indexed stochastic processes $X_{t,T}$ ($t = 0, \ldots, T-1, T > 0$) with mean zero is in the class of locally stationary wavelet processes (LSW processes) if there exists a representation

$$(2.1) \qquad X_{t,T} = \sum_{j=-\infty}^{-1} \sum_{k=0}^{T-1} w_{jk;T} \psi_{jk}(t) \xi_{jk},$$

where the infinite sum is to be understood in the mean-square sense, $\{\psi_{jk}(t) = \psi_{j,k-t}\}_{jk}$ with $j < 0$ is a discrete nondecimated family of wavelets based on a mother wavelet $\psi(t)$ of compact support, and such that the following conditions are satisfied:

1. $\xi_{jk}$ is a random orthonormal increment sequence such that $\mathrm{E}\xi_{jk} = 0$ and $\mathrm{Cov}(\xi_{jk}, \xi_{\ell m}) = \delta_{j\ell}\delta_{km}$ for all $j, \ell, k, m$, where $\delta_{j\ell} = 1$ if $j = \ell$ and 0 elsewhere.
2. For each $j \leq -1$, there exists a function $W_j(z)$ on $(0,1)$ possessing the following properties:
   (a) $\sum_{j=-\infty}^{-1} |W_j(z)|^2 \leq \overline{C} < \infty$ uniformly in $z \in (0,1)$;
   (b) there exists a sequence of constants $C_j$ such that for each $T$

$$(2.2) \qquad \sup_{k=0,\ldots,T-1} \left| w_{jk;T} - W_j\left(\frac{k}{T}\right) \right| \leq \frac{C_j}{T};$$

   (c) the total variation of $W_j^2(z)$ is bounded by $L_j$, that is,

$$\mathrm{TV}(W_j^2) := \sup\left\{ \sum_{i=1}^{I} |W_j^2(a_i) - W_j^2(a_{i-1})| : 0 < a_0 < \cdots < a_I < 1, \right.$$

$$(2.3) \qquad \left. I \in \mathbb{N} \right\}$$

$$\leq L_j,$$

   (d) the constants $C_j$ and $L_j$ are such that

$$(2.4) \qquad \sum_{j=-\infty}^{-1} \mathcal{L}_j(\mathcal{L}_j L_j + C_j) \leq \rho < \infty,$$

   where $\mathcal{L}_j = |\mathrm{supp}\,\psi_{j0}| = (2^{-j} - 1)(\mathcal{L}_{-1} - 1) + 1$.

LSW processes use wavelets to decompose a stochastic process with respect to an orthogonal increment process in the time-scale plane. Due to the overcompleteness of the nondecimated system, a given LSW processes does not determine the sequence $\{w_{jk;T}\}$ uniquely. However, we can build a theory which ensures the existence of a unique wavelet spectrum (in a sense defined after Proposition 1 below). This property is a consequence of the local stationarity setting which introduces a *rescaled time* $z = t/T \in (0,1)$ on which $W_j(z)$ is defined. The rescaled time permits increasing amounts of data about the local structure of $W_j(z)$ to be collected as the observed time $T$ tends to infinity. Even though a given LSW process does not determine the sequence $\{w_{jk;T}\}$ uniquely, the model allows the (asymptotic) identification of the model coefficients determined by uniquely defined $W_j^2(z)$. Then, the *evolutionary wavelet spectrum* of an LSW process $\{X_{t,T}\}_{t=0,\ldots,T-1}$, with respect to $\psi$, is given by

$$(2.5) \qquad S_j(z) = |W_j(z)|^2, \qquad z \in (0,1),$$

and is such that, by definition of the process, $\sum_{j=-\infty}^{-1} S_j(z) < \infty$ uniformly in $z \in (0,1)$.

The evolutionary wavelet spectrum $S_j(z)$ is related to the time-dependent autocorrelation function of the LSW process. Observe that the autocovariance function of an LSW process can be written as

$$c_{X,T}(z,\tau) = \mathrm{Cov}(X_{[zT],T}, X_{[zT]+\tau,T})$$

for $z \in (0,1)$ and $\tau$ in $\mathbb{Z}$, and where $[\cdot]$ denotes the integer part of a real number. The next result shows that this autocovariance converges asymptotically to a *local autocovariance* defined by

$$(2.6) \qquad c_X(z,\tau) = \sum_{j=-\infty}^{-1} S_j(z)\Psi_j(\tau),$$

where $\Psi_j(\tau) = \sum_{k=-\infty}^{\infty} \psi_{jk}(0)\psi_{jk}(\tau)$ is the *autocorrelation wavelet* function.

PROPOSITION 1.  *Under the assumptions of Definition 1, if $T \to \infty$, then*

$$\sum_{\tau=-\infty}^{\infty} \int_0^1 dz |c_{X,T}(z,\tau) - c_X(z,\tau)| = O(T^{-1})$$

*for all LSW process.*

Appendix A presents some properties of the autocorrelation wavelet system appearing in (2.6). Like wavelets themselves, this system enjoys good localization properties. Consequently, we observe that equation (2.6) is a multiscale decomposition of the autocovariance structure of the process over time: the larger the wavelet spectrum $S_j(z)$ is at a particular scale $j$ and point $z$ in the rescaled time, the more dominant is the contribution of scale $j$ in the variance at time $z$. Thus, the evolutionary wavelet spectrum describes the distribution of the (co)variance at a particular scale and time location.

Moreover, we recall in Appendix A that $\{\Psi_j\}$ is a linearly independent system. Therefore, since the autocovariance function converges to the local autocovariance in the sense of Proposition 1, the coefficients $S_j(z)$ in (2.6) are asymptotically the unique wavelet representation of the second order structure of the time series.

It is worth mentioning that a stationary process with an absolutely summable autocovariance function is an LSW process [Nason et al. (2000), Proposition 3]. Stationarity is characterized by a wavelet spectrum which is constant over time: $S_j(z) = S_j$ for all $z \in (0,1)$. However, our motivation for studying LSW processes lies in the modeling of time-varying spectra. The regularity of the wavelet spectrum in time is determined by the smoothness of $W_j(z)$ with respect to $z$. In Nason et al. (2000), this function is assumed to be Lipschitz continuous in time. In our definition of LSW processes, we only require the total variation of $W_j^2$ to be bounded. This weaker assumption is considered not only in order to work with less strict assumptions, but also to allow a discontinuous evolution of the wavelet spectrum in time. Figure 1 shows a simulated example of such a nonstationary process.

## 3. A first estimator of the wavelet spectrum.

3.1. *The corrected wavelet periodogram.* An estimator of the wavelet spectrum is constructed by taking the squared empirical coefficients from the nondecimated transform:

$$I_{j;T}\left(\frac{k}{T}\right) = \left(\sum_{t=0}^{T-1} X_{t,T} \psi_{jk}(t)\right)^2, \qquad j = -1, \ldots, -\log_2 T; k = 0, \ldots, T-1.$$

$I_{j;T}(z)$ is called the *wavelet periodogram*, as it is analogous to the formula for the classical periodogram in traditional Fourier spectral analysis of stationary processes [Brillinger (1975)].

Some asymptotic properties of this estimator have been studied by Nason et al. (2000), who showed that the wavelet periodogram is *not* an asymptoticaly unbiased estimator of the wavelet spectrum. Indeed, Proposition 4 of Nason et al. (2000) states that for all fixed scales $j < 0$,

(3.1) $$\mathrm{E}\, I_{\ell;T}(z) - \sum_{\ell=-\log_2 T}^{-1} A_{j\ell} S_\ell(z) = O(T^{-1}),$$

uniformly in $z \in (0,1)$, where the matrix $A = (A_{j\ell})_{j,\ell<0}$ is defined by

$$A_{j\ell} := \langle \Psi_j, \Psi_\ell \rangle = \sum_\tau \Psi_j(\tau) \Psi_\ell(\tau).$$

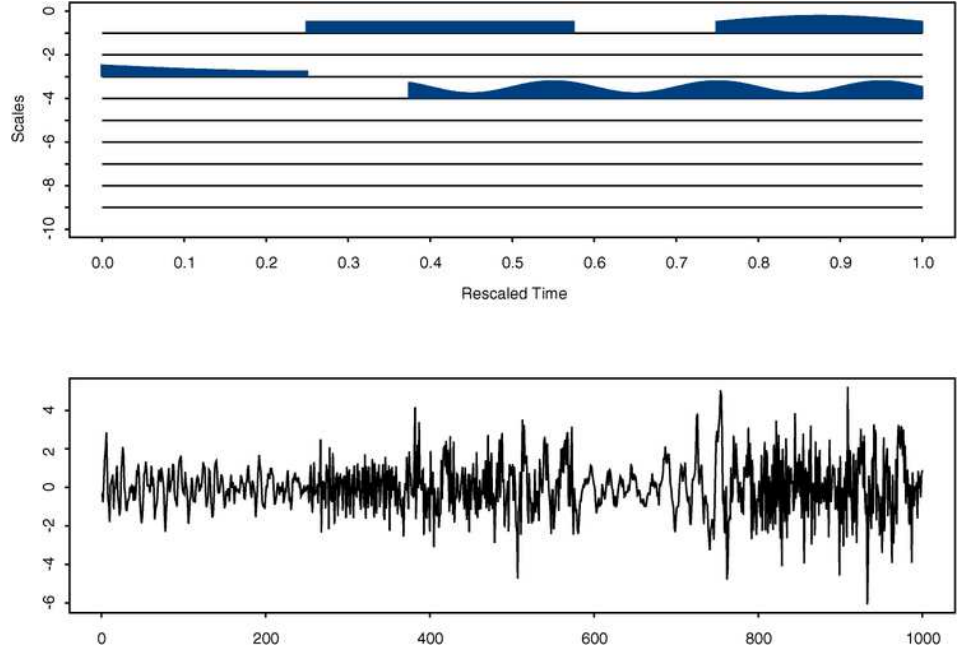

Fig. 1. *The upper figure is an example of theoretical spectrum $S_j(z)$. This spectrum is used in the lower figure to simulate a locally stationary wavelet process of length $T = 1000$. This simulation uses Gaussian innovations $\xi_{jk}$ and nondecimated Haar wavelets.*

Note that the matrix $A_{j\ell}$ is not simply diagonal since the autocorrelation wavelet system $\{\Psi_j\}$ is not orthogonal. Nason et al. (2000) proved the invertibility of $A$ if $\{\Psi_j\}$ is constructed using Haar wavelets. If other compactly supported wavelets are used, numerical results suggest that the invertibility of $A$ still holds, but a complete proof of this result has not yet been established. As we need the invertibility of $A$ in results which follow, we hereafter restrict ourselves to Haar wavelets, but conjecture that all results remain valid for more general Daubechies wavelets [Daubechies (1992)].

Equation (3.1) motivates the definition of a *corrected wavelet periodogram*,

$$(3.2) \qquad L_{j;T}\left(\frac{k}{T}\right) = \sum_{\ell=-\log_2 T}^{-1} (A_T)^{-1}_{j\ell}\left(\sum_{t=0}^{T-1} X_{t,T}\psi_{\ell k}(t)\right)^2,$$

where $A_T = (A_{j\ell})_{-\log_2 T \leq j,\ell \leq -1}$. The corrected wavelet periodogram $L_{j;T}$ is a preliminary tool for constructing an asymptotically consistent estimator of the evolutionary wavelet spectrum. To this end, it needs to be smoothed in time. This question is addressed in the following.

REMARK 1. The asymptotic bias of the wavelet periodogram is a consequence of the overcompleteness of the nondecimated wavelet system $\{\psi_{jk}\}$.

One could ask if it would not be easier to define LSW processes using a decimated wavelet system because, for this system, the matrix $A$ reduces to the identity. Unfortunately, the answer is negative: the use of nondecimated wavelets, as described in von Sachs et al. (1997), would not allow the local autocovariance function to be written as a wavelet-type transform of an evolutionary spectrum, as in (2.6). Moreover, classical stationary processes are not included in the model based on decimated wavelets.

3.2. *The preliminary estimator and its properties.* Suppose we want to estimate $S_j(z_0)$ from observations $\underline{X}_T = (X_{0,T}, \ldots, X_{T-1,T})'$. The estimator studied below takes the following form:

$$(3.3) \quad Q_{j,\mathcal{R};T} = |\mathcal{R}T|^{-1} \sum_{k \in \mathcal{R}T} \left\{ L_{j;T}\left(\frac{k}{T}\right) + z_{j,k;T} \right\}, \qquad j = -1, -2, \ldots,$$

where $z_{j,k;T}$ are i.i.d. Gaussian random variables of mean zero and variance $C^2 2^j$, independent from $\underline{X}_T$ for a given constant $C^2$, $\mathcal{R}$ is an interval in $(0,1)$ that contains the point $z_0$ and $k \in \mathcal{R}T$ means that $k/T \in \mathcal{R}$. The estimator (3.3) is essentially the average of the corrected wavelet periodogram over the interval $\mathcal{R}$. The reason for adding a "noise process" $z_{j,k;T}$ in our estimator is for the sake of regularization, since the process $\underline{X}_T$ is not guaranteed to be invertible. In other words, the presence of the additive Gaussian variable in the estimator $Q_{j,\mathcal{R};T}$ allows consistent estimation of more general processes for which the wavelet spectrum $S_j(z)$ is not bounded away from zero. Note that this regularization technique does not add any systematic bias to the resulting estimator since in (3.3), an average is taken over the zero-mean Gaussian variables $z_{j,k;T}$. That procedure is analogous to the regularization techniques for ill-posed inverse problems such as, for instance, in ridge regression or Tikhonov regularization; see also Neumann (1996) for a similar technique in the context of stationary time series.

Of course, the choice of the interval $\mathcal{R}$ around $z_0$ is crucial in this estimation. This question will be addressed in the next section. First, we derive some useful properties of $Q_{j,\mathcal{R};T}$ as an estimator of

$$(3.4) \qquad Q_{j,\mathcal{R}} = |\mathcal{R}|^{-1} \int_{\mathcal{R}} dz \, S_j(z).$$

The statistical properties of $Q_{j,\mathcal{R};T}$ are now derived under a set of assumptions.

ASSUMPTION 1. The autocovariance function $c_{X,T}$ and the local autocovariance function $c_X$ of the LSW process are such that

$$(3.5) \qquad \|c_{X,T}\|_{1,\infty} := \sum_{\tau=-\infty}^{\infty} \sup_{t=0,\ldots,T-1} \left| c_{X,T}\left(\frac{t}{T}, \tau\right) \right|$$

is bounded independently of $T$

and

(3.6) $$\|c_X\|_{1,\infty} := \sum_{\tau=-\infty}^{\infty} \sup_{z \in (0,1)} |c_X(z,\tau)| < \infty.$$

This assumption is needed to control the spectral norm of the covariance matrix of the process (Lemma B.3 in Appendix B). For a stationary process, it reduces to absolute summability of the autocovariance of the process (short memory property).

ASSUMPTION 2. *There exists an $\varepsilon > 0$ such that for all $z \in (0,1)$, $\sum_{j=-\infty}^{-1} S_j(z) \geq \varepsilon$.*

According to equation (2.6), the sum over scales of $S_j(z)$ is the local variance of the process at time $[zT]$ and this assumption states that the local variance of the process is bounded away from zero.

ASSUMPTION 3. *The increment process $\{\xi_{jk}\}$ in Definition 1 is Gaussian.*

This assumption allows substantial simplifications in the proofs. It is also assumed to establish some results in Nason et al. (2000) and Fryźlewicz et al. (2003).

ASSUMPTION 4. *The evolutionary wavelet spectrum $S_j(z)$ is such that*

$$\sum_{\ell=-\infty}^{-\log_2(T)-1} \sup_{z \in (0,1)} S_\ell(z) = O(T^{-1}).$$

In the definition of the corrected wavelet periodogram (3.2), all scales $0 > j > -\infty$ are implicitly included due to the definition of $X_{t,T}$. The last assumption is used in order to control the remainder of the estimation bias at all scales lower than $-\log_2 T$.

The following proposition describes the asymptotic properties of $Q_{j,\mathcal{R};T}$.

PROPOSITION 2. *Suppose Assumptions 1–4 hold true. For all LSW processes (Definition 1) and all $\mathcal{R} \subseteq (0,1)$,*

(3.7) $$\mathrm{E}\, Q_{j,\mathcal{R};T} - Q_{j,\mathcal{R}} = \frac{K_0 2^{j/2}\sqrt{T}}{|\mathcal{R}T|} \sum_{m=-\log_2 T}^{-1} \mathcal{L}_m \,\mathrm{TV}(S_m) + O(2^{j/2}|\mathcal{R}T|^{-1})$$

$$= O\!\left(\frac{2^{j/2}}{\sqrt{T}}\right)$$

for all $j = -1, \ldots, -J_T$ with $J_T = O(\log_2 T)$ and where $K_0$ is a constant independent of $j, T$ and $|\mathcal{R}|$. Moreover, under Assumptions 1–4, the variance $\sigma^2_{j,\mathcal{R};T} = \operatorname{Var} Q_{j,\mathcal{R};T}$ is such that

$$\frac{C^2 2^j}{|\mathcal{R}T|} \leq \sigma^2_{j,\mathcal{R};T} \leq \left(C^2 + \frac{c^2}{|\mathcal{R}|}\right) \frac{2^j}{|\mathcal{R}T|}$$

for all $T$, for all $j = -1, \ldots, -J_T$ with $J_T = o_T(\log_2 T)$ and where $c^2 = 2K_2^2 \|c_X\|^2_{1,\infty}$ with $K_2$ a constant that depends only on the wavelet $\psi$.

The proof of this proposition is in Appendix B.3. Note that the squared bias and the variance of the estimator have the same rate of convergence. This phenomenon is due to the nonstationary behavior of the process. Indeed, for a stationary process, the total variation of $S_m$ is zero at all scales and the rate of the bias is then $T^{-1}$. This is not the case for a general nonstationary process: when the wavelet spectrum is not constant over time, an additional term resulting from nonstationarity considerably reduces this rate of convergence. Moreover, even if we are dealing with a *local* estimator of the wavelet spectrum at a fixed scale $j < 0$ and a fixed time interval $\mathcal{R}$, the nonstationarity term in the bias involves the variation of the *global* wavelet spectrum. This may be observed in equation (3.7), which involves a sum over all scales $m = -1, \ldots, -\log_2 T$ and the total variation of all $S_m$ over the whole rescaled time interval $(0, 1)$.

This slow rate of convergence of the bias poses a problem for the establishment the asymptotic normality of $Q_{j,\mathcal{R};T}$. In the next proposition, we circumvent this problem and derive a nonasymptotic exponential bound for the deviation of $Q_{j,\mathcal{R};T}$.

PROPOSITION 3. *Assume that Assumptions 1–4 hold. If $\sigma^2_{j,\mathcal{R},T} = \operatorname{Var} Q_{j,\mathcal{R};T}$, then for all $\eta > 0$ and all scales $j = -1, \ldots, -J_T$, where $J_T = O(\log_2 T)$,*

$$\Pr(|Q_{j,\mathcal{R};T} - Q_{j,\mathcal{R}}| > 2\sigma_{j,\mathcal{R},T}\eta)$$
$$\leq c_0 \exp\left\{-\frac{1}{16} \cdot \eta^2 \Big/ \left[1 + \frac{2\eta L_j}{|\mathcal{R}T|\sigma_{j,\mathcal{R},T}} + \frac{2^{j/2}\eta(K_2\|c_X\|_{1,\infty} + K_3)}{|\mathcal{R}|\sqrt{T}\sigma_{j,\mathcal{R},T}}\right]\right\}$$

*with the positive constants $c_0 = 3 + e$, $K_2$ as in Proposition 2 and $K_3$ depending on the wavelet $\psi$ and the constants $\rho, \overline{C}$ given in Definition 1.*

The proof of this proposition appears in Appendix B.4. This proposition gives a nonasymptotic approximation for the deviation of $Q_{j,\mathcal{R};T}$. This result is exploited in the next section in order to choose the interval $\mathcal{R}$ in an adaptive way. From an asymptotic viewpoint, that is, as $T \to \infty$, we note that this exponential bound does not tend to zero, meaning that the standardized statistic $Q_{j,\mathcal{R},T}$ is asymptotically nondegenerate. This phenomenon is well known in the context of pointwise estimation; see Lepski (1990) and Brown and Low (1996). In order to have a consistent result when $T \to \infty$, it is then necessary to require that $\eta = \eta_T$ grows with $T$. The appropriate rate for $\eta_T$ is derived in the next corollary. The proof is given in Appendix B.4 and is essentially based on the bounds derived in Proposition 2.

COROLLARY 1. *Under the assumptions of Propositions 2 and 3, if $k_T$ tends to infinity and is such that $J_T \cdot \exp(-k_T) = o_T(1)$, then there exists a $T_0 > 1$ such that for all $T \geq T_0$,*

$$\Pr\left(\sup_{-J_T \leq j < 0} |Q_{j,\mathcal{R};T} - Q_{j,\mathcal{R}}| \geq k_T \sqrt{(1 + c^2/|\mathcal{R}|)/|\mathcal{R}T|}\right) = o_T(1),$$

*where $c^2$ is as in the assertion of Proposition 2.*

REMARK 2. An example of admissible rates is $J_T \sim \log_2 T$ and $k_T \sim \log_2 T$. The sequence $k_T$ will play a crucial role in Section 4.

REMARK 3. The results are proved under the assumption that the increments considered in the definition of LSW processes are Gaussian (Assumption 3). This assumption allows substantial simplifications in the proofs. For practical applications, we believe that this assumption is not unrealistic and the class of Gaussian LSW processes is rich enough, as can be observed from the wide range of applications that are treated in Nason et al. (2000), Fryźlewicz et al. (2003), Oh et al. (2003), Woyte et al. (2007) and Van Bellegem and von Sachs (2004), for instance. However, it still seems interesting to see how the above results can be extended to the non-Gaussian case. A careful reading of the proof of Proposition 3 shows that the crucial point is to establish an exponential inequality for quadratic forms of the increments. In our proof of Proposition 3, we use the inequality established by Dahlhaus and Polonik (2002) on the quadratic form of Gaussian random variables. Other exponential inequalities have been established for non-Gaussian random variables; see, for instance, Dahlhaus (1988) or Spokoiny (2001, 2002). Another example of an exponential inequality for dependant data is derived in van de Geer (2002).

3.3. *Estimation of the variance.* The main drawback of Proposition 3 is that the deviation result depends on the variance $\sigma^2_{j,\mathcal{R},T} = \text{Var}\, Q_{j,\mathcal{R};T}$ which is typically unknown. The goal of the following derivation is to propose a preliminary estimator $\tilde{\sigma}^2_{j,\mathcal{R},T}$ of $\sigma^2_{j,\mathcal{R},T}$ such that Proposition 3 can still be used with $\tilde{\sigma}^2_{j,\mathcal{R},T}$.

The variance $\sigma^2_{j,\mathcal{R},T}$ depends on the unknown autocovariance function of the LSW process in the following way [see Lemma B.1 with equation (B.9)]:

$$\sigma^2_{j,\mathcal{R},T} = 2\|U'_{j,\mathcal{R};T}\Sigma_T\|^2_2 + \frac{C^2 2^j}{|\mathcal{R}T|},$$

where $\Sigma_T$ is the $T \times T$ (non-Toeplitz) covariance matrix of the LSW process $(X_{0,T}, \ldots, X_{T-1,T})'$, and $U_{j,\mathcal{R};T}$ is the $T \times T$ matrix with entry $(s,t)$ equal to

$$U^{(j)}_{st} = |\mathcal{R}T|^{-1} \sum_{\ell=-\log_2 T}^{-1} A^{-1}_{j\ell} \sum_{k \in \mathcal{R}T} \psi_{\ell k}(s)\psi_{\ell k}(t).$$

We also denote by $\sigma_{s,s+u}$ the entry $(s, s+u)$ of the matrix $\Sigma_T$.

We will estimate $\sigma^2_{j,\mathcal{R},T}$ by

$$\tilde{\sigma}^2_{j,\mathcal{R},T} = 2\|U'_{j,\mathcal{R};T}\tilde{\Sigma}_T\|^2_2 + \frac{C^2 2^j}{|\mathcal{R}T|},$$

where $\tilde{\Sigma}_T$ is an estimate of the covariance matrix $\Sigma_T$. A first idea is to define the elements $\tilde{\sigma}_{s,s+u}$ of $\tilde{\Sigma}_T$ by plugging $Q_{j,\mathcal{R};T}$ into the local autocovariance function (2.6), that is,

$$\tilde{\sigma}_{s,s+u} = \sum_{j=-\log_2 T}^{-1} Q_{j,\mathcal{R}(s);T}\Psi_j(u),$$

where $\mathcal{R}(s)$ denotes an interval which contains the time point $s/T$. However, the convergence in probability of $\tilde{\sigma}_{s,s+u}$ to $\sigma_{s,s+u}$ is not faster than the rate of $\sigma_{s,s+u}$ itself and we need to modify the estimator in the following two ways.

(i) Assumption 1 indicates that the covariance $|\sigma_{s,s+u}|$ is small for large $|u|$. We set $\tilde{\sigma}_{s,s+u}$ to zero when $|u| \geq M_T$ for an appropriate sequence $M_T$ tending to infinity with $T$.
(ii) It is necessary to control the distance in rescaled time between the spectrum $S_j(z)$, for $z \in \mathcal{R}(s)$, and $S_j(s/T)$. To do so, we allow the window $\mathcal{R}(s)$ to depend on $T$, which is denoted by $\mathcal{R}_T(s)$, in such a way that its length $|\mathcal{R}_T|$ shrinks to zero when $T$ tends to infinity. This is analogous to the estimation of a regression function by kernel smoothing, where the window usually depends on the length of the data set.

With these two ingredients, we propose to estimate $\sigma_{s,s+u}$ by

(3.8) $$\tilde{\sigma}_{s,s+u} = \sum_{j=-\log_2 T}^{-1} Q_{j,\mathcal{R}_T(s);T}\Psi_j(u)\mathbb{I}_{|u|\leq M_T}$$

and the following assumption makes precise the appropriate rates for the sequences $|\mathcal{R}_T|$ and $M_T$.

ASSUMPTION 5. The sequence $J_T$ is such that $J_T = o_T(\log_2 T)$. The length of $\mathcal{R}_T$ tends to zero such that $2^{J_T}|\mathcal{R}_T| = o_T(1)$. The sequence $k_T$ (which appears in Corollary 1) tends to infinity such that $J_T \exp(-k_T\sqrt{|\mathcal{R}_T|}) = o_T(1)$. Finally, the sequence $M_T$ [involved in the preliminary estimator for the variance—see (3.8)] tends to infinity such that

$$2^{J_T}|\mathcal{R}_T|^{-1}T^{-1/2}M_T k_T \log_2^3 T = o_T(1).$$

Admissible rates for this last assumption are, for example, $J_T \sim \log_2 \log_2^2 T$, $k_T \sim \log^2 T$, $|\mathcal{R}_T| \sim \log_2^{-3} T$ and $M_T \sim \log_2^\alpha T$ with $\alpha > 0$. It is worth mentioning that with this assumption, $|\mathcal{R}_T|$ shrinks to zero in the *rescaled* time, whereas in the *observed* time, the interval length $|T\mathcal{R}_T|$ tends to infinity. This means that our estimate of $S_j(s/T)$ is built using an increasing amount of data in the observed time, but, at the same time, with an average around $S_j(s/T)$ in the rescaled time on a shrinking segment around $s/T$.

The next proposition shows that on the random set where the estimator $Q_{j,\mathcal{R}_T(s);T}$ is near $Q_{j,\mathcal{R}_T(s)}$, the estimator (3.8) has a good quality. Our proof of this proposition may be found in Appendix B.5 and needs the following technical assumption, which is a slightly stronger condition than point 2(a) of Definition 1, in the sense that we need to control the decay of $S_j(z)$ with respect to $j$ and uniformly in $z$.

ASSUMPTION 6. The local autocovariance function $c(z,\tau)$ is such that

$$\sum_{u=-\infty}^{\infty} \sup_z |c_X(z,u)| \mathbb{I}_{|u|>M_T} = o_T(2^{-J_T}).$$

This last assumption on the decay of the local autocovariance function uniformly in $z$ is very sensible in the context of short-memory stationary processes [in that case, $c(z,u)$ does not depend on $z$]. With the rates specified above, a typical condition is to assume $|c_X(z,u)| \leq c \cdot r^{|u|}$ uniformly in $z \in (0,1)$ with $0 \leq r < 1$.

PROPOSITION 4. *Suppose Assumptions 1–6 hold. There then exists a positive number $T_0$ and a random set $\mathcal{A}$ independent of $j$ and such that $\Pr(\mathcal{A}) \geq 1 - o_T(1)$ and*

$$|Q_{j,\mathcal{R}_T(s);T} - Q_{j,\mathcal{R}_T(s)}| \leq K_2 \|c_X\|_{1,\infty} k_T \sqrt{T} |\mathcal{R}_T T|^{-1}$$

*for all $T > T_0$. Moreover, on $\mathcal{A}$,*

(3.9) $$2^{J_T-j}T|\tilde{\sigma}^2_{j,\mathcal{R},T} - \sigma^2_{j,\mathcal{R},T}| = o_P(1)$$

*holds for all $j = -1,\ldots,-J_T$, where $o_P(1)$ does not depend on $\mathcal{R}$.*

Finally, Proposition 4 together with Proposition 3 leads to the following result, which will be used to construct the pointwise adaptive estimator in Section 4.

THEOREM 1. *Suppose Assumptions* 1–6 *hold. There then exists a* $\gamma_T = o_T(1)$ *and a positive number* $T_0$ *such that for all* $T > T_0$,

$$\Pr(|Q_{j,\mathcal{R};T} - Q_{j,\mathcal{R}}| > 2\tilde{\sigma}_{j,\mathcal{R},T}\eta')$$
$$\leq c_0 \exp\left\{-\frac{1}{16} \cdot \eta^2 \Big/ \left[1 + \frac{2\eta L_j}{|\mathcal{R}T|\sigma_{j,\mathcal{R},T}} + \frac{2^{j/2}\eta(K_2\|c_X\|_{1,\infty} + K_3)}{|\mathcal{R}|\sqrt{T}\sigma_{j,\mathcal{R},T}}\right]\right\} + o_T(1)$$

*for all* $j = -1,\ldots,-J_T$, *where* $\eta' = \eta\sqrt{1-\gamma_T}$ *and the positive constants* $c_0, K_2, K_3$ *are defined as in the assertion of Propositions* 2 *and* 3.

REMARK 4. Theorem 1 gives an approximation of the distribution of the normalized loss $|Q_{j,\mathcal{R};T} - Q_{j,\mathcal{R}}|/\tilde{\sigma}_{j,\mathcal{R},T}$. This depends on the unknown quantities $\|c_X\|_{1,\infty}$ and $\rho$ [cf. (2.4)]. These two quantities may be understood as nuisance parameters of the problem, depending on the global spectrum. The estimation of these quantities is based on a preliminary smoothing of $L_{j;T}(z)$ with respect to $z$, which we denote by $L^*_{j;T}(z)$. Here, we think about using a kernel smoothing procedure, or a wavelet transform shrinkage as studied in Nason et al. (2000). A preliminary estimate of $\|c_X\|_{1,\infty}$ is then obtained by plugging $L^*_{j;T}(z)$ into $\|c_X\|_{1,\infty}$, [cf. (2.6) and (3.6)]. Next, the preliminary estimation of $\rho$ necessitates the estimation of $\text{TV}(S_j)$ [cf. (2.3)]. We estimate $\text{TV}(S_j)$ by $\sum_i |L^*_{j;T}(z_i^{\max}) - L^*_{j;T}(z_i^{\min})| + |L^*_{j;T}(z_i^{\max}) - L^*_{j;T}(z_{i+1}^{\min})|$, where the sum is over the local minima and maxima of $L^*_{j;T}(z)$, with $z_i^{\max} < z_{i+1}^{\min} < z_{i+1}^{\max}$ for all $i$.

REMARK 5. The estimator (3.3) also involves a constant $C^2$. In view of Proposition 2 on the variance of the estimator, that constant should ideally be close to $c^2 = 2K_2^2\|c_X\|_{1,\infty}$. Because $\|c_X\|_{1,\infty}$ is unknown, it is estimated in practice by $\sum_s \sup_u \tilde{\sigma}_{s,s+u}$.

**4. Pointwise adaptive estimation.** The question of how to choose the best segment $\mathcal{R}$ in the estimator (3.3) arises and the goal of this section is to provide a data-driven procedure to select $\mathcal{R}$ automatically.

The proposed method goes back to the pointwise adaptive estimation theory of Lepski (1990); see also Lepski and Spokoiny (1997) and Spokoiny (1998). Suppose that the wavelet spectrum $S_j(z_0)$ is well approximated by the averaged spectrum $Q_{j,\mathcal{U}}$ for a given interval $\mathcal{U}$ containing the reference point $z_0$. The idea of the procedure is to consider a second interval $\mathcal{R}$ containing $\mathcal{U}$ and to test whether $Q_{j,\mathcal{R}}$ differs significantly from $Q_{j,\mathcal{U}}$. As we describe below, this test procedure is based on Proposition 3 or Theorem 1. If there exists a subset $\mathcal{U}$ of $\mathcal{R}$ such that $|Q_{j,\mathcal{R}} - Q_{j,\mathcal{U}}|$ is significantly different from zero, then we reject the hypothesis of homogeneity of the wavelet spectrum $S_j(z)$ on $z \in \mathcal{R}$. Finally, the adaptive estimator corresponds to the largest interval $\mathcal{R}$ such that the hypothesis of homogeneity is not rejected.

This section contains a formal description of this algorithm and derives some properties of the estimator.

4.1. *Testing homogeneity.* Let $\mathcal{R}$ be an interval containing $z_0$, $\mathcal{U}$ a subset of $\mathcal{R}$ and define

$$(4.1) \qquad \Delta_j(\mathcal{R},\mathcal{U}) = |Q_{j,\mathcal{R}} - Q_{j,\mathcal{U}}|.$$

Under Assumptions 1–4, Proposition 3 implies that

$$\Pr[|Q_{j,\mathcal{R},T} - Q_{j,\mathcal{U},T}| > \Delta_j(\mathcal{R},\mathcal{U}) + 2\eta(\sigma_{j,\mathcal{R},T} + \sigma_{j,\mathcal{U},T})k_T]$$
$$\leq h(\mathcal{U},\eta) + h(\mathcal{R},\eta)$$

with

$$h(\mathcal{R},\eta) = c_0 \exp\left\{-\frac{1}{16} \cdot (\eta^2 k_T^2) \Big/ \left[1 + \frac{2\eta k_T}{|\mathcal{R}T|\sigma_{j,\mathcal{R},T}} L_j + \frac{2^{j/2}\eta k_T}{|\mathcal{R}|\sqrt{T}\sigma_{j,\mathcal{R},T}}(K_2 \|c_X\|_{1,\infty} + K_3)\right]\right\}$$

and where the sequence $k_T$ is such that $J_T \cdot \exp(-k_T) = o_T(1)$ (see Corollary 1). Under the assumption that the wavelet spectrum $S_j$ is homogeneous on the segment $\mathcal{R}$, the difference $\Delta_j(\mathcal{R},\mathcal{U})$ is negligible. Then, as a test rule, we reject the homogeneity hypothesis on $\mathcal{R}$ if there exists a subset $\mathcal{U} \subset \mathcal{R}$ such that $|Q_{j,\mathcal{R};T} - Q_{j,\mathcal{U};T}| > 2\eta(\sigma_{j,\mathcal{R},T} + \sigma_{j,\mathcal{U},T})k_T$ for a given $\eta$.

In the case where the variances $\sigma_{j,\mathcal{R},T}$ and $\sigma_{j,\mathcal{U},T}$ are unknown, they may be estimated as in Section 3.3 above.

In practice, we choose a set $\Lambda$ of interval candidates $\mathcal{R}$. Then, for each candidate $\mathcal{R}$, we apply the homogeneity test with respect to a given set $\wp(\mathcal{R})$ of subintervals $\mathcal{U}$ of $\mathcal{R}$.

ASSUMPTION 7. In the estimation procedure described below, we assume the following properties on the test sets $\Lambda$ and $\wp(\mathcal{R})$:

1. for all $\mathcal{R}$, the shortest interval of $\wp(\mathcal{R})$ is of length at least $\delta > 0$;
2. the cardinality of $\wp(\mathcal{R})$ is such that $\sharp(\wp(\mathcal{R})) \leq |\mathcal{R}T|^{(\alpha\sqrt{\delta K_1})/(K_2\|c_X\|_{1,\infty}+K_3)}$ for some $0 < \alpha < 1$;

3. when we test the homogeneity of the wavelet spectrum on $\mathcal{R}$, we assume that there exists a subinterval $\mathcal{U} \in \wp(\mathcal{R})$ such that $\mathcal{U} \subset \mathcal{R}$ and $\mathcal{U}$ contains $z_0$.

REMARK 6 (Test sets). We give an example of sets $\Lambda$ and $\wp(\mathcal{R})$. For each scale $j < 0$, the corrected wavelet spectrum (3.2) is evaluated on a grid $k/T$, $r = 0, \ldots, T-1$ in time. We can then choose the set $\Lambda$ as

$$\Lambda = \{[r_0/T, r_1/T] : r_0 < [z_0 T] < r_1\}$$

for $r_0, r_1 \in \{0, T-1\}$. Nevertheless, in order to reduce the computational effort, we shrink the cardinality of $\Lambda$ following the method of Spokoiny (1998). More precisely, we first select two sets $\mathcal{K}_m = \{r_m : r_m \leq [z_0 T]\}$ and $\mathcal{K}_n = \{r_n : r_n \geq [z_0 T]\}$ which both contain less than $T$ points and then set

$$\Lambda = \{[r_m/T, r_n/T] : r_m \in \mathcal{K}_m, r_n \in \mathcal{K}_n\}.$$

Then, one possibility for defining $\wp(\mathcal{R})$ is to consider

$$\wp(\mathcal{R}) = \{[r_-/T, r_+/T] : r_-, r_+ \in \mathcal{K}_m \cup \mathcal{K}_n\}.$$

We refer to Spokoiny (1998) for the details of this construction.

4.2. *The estimation procedure.* The estimation procedure simply starts with the smallest interval in $\Lambda$, assuming that the wavelet spectrum is homogeneous on this short interval. It then iteratively selects longer intervals in $\Lambda$ until the homonegeneity assumption is rejected. Finally, the adaptive segment $\tilde{\mathcal{R}}$ is the longest segment $\mathcal{R}$ of $\Lambda$ for which the homogeneity test is not rejected,

$$\tilde{\mathcal{R}} = \arg\max_{\mathcal{R} \in \Lambda}\{|\mathcal{R}| \text{ such that } |Q_{j,\mathcal{R};T} - Q_{j,\mathcal{U};T}| \leq 2\eta(\sigma_{j,\mathcal{R},T} + \sigma_{j,\mathcal{U},T})k_T$$

(4.2)
$$\text{for all } \mathcal{U} \subset \wp(\mathcal{R})\}.$$

The adaptive estimator of $S_j(z_0)$ is then defined by

$$\tilde{S}_j(z_0) = Q_{j,\tilde{\mathcal{R}},T}. \tag{4.3}$$

In the case where the variances $\sigma_{j,\mathcal{R},T}$ and $\sigma_{j,\mathcal{U},T}$ are unknown, they may be estimated as in Section 3.3 above. In that case, the homogeneity test is based on Theorem 1 and the modification of the following results is straightforward. The proofs are longer, however, but the technique in the proof of Theorem 1 can be used to transfer the problem with estimated variances to the problem with known variances $\sigma_{j,\mathcal{R},T}$ and $\sigma_{j,\mathcal{U},T}$.

4.3. *Properties of the estimator in homogeneous regions.* The next result quantifies the $\ell_p$-risk ($p \geq 2$) when the wavelet spectrum $S_j(z)$ is homogeneous on $z \in \mathcal{R}$. To define this concept of homogeneity, we introduce the bias

$$b(\mathcal{R}) := \sup_{z \in \mathcal{R}} |S_j(z) - Q_{j,\mathcal{R}}|,$$

which measures how well the wavelet spectrum $S_j$ is approximated by $Q_{j,\mathcal{R}}$ on $z \in \mathcal{R}$. We say that the spectrum is *homogeneous* (or *regular*) on $\mathcal{R}$ if the inequality

(4.4) $$b(\mathcal{R}) \leq C_j \sigma_{j,\mathcal{R},T} k_T$$

holds with

(4.5) $$C_j = 2^{-j/2} \sqrt{\alpha + p}$$

for a positive real constant $\alpha$. In the inequality (4.4), $\sigma_{j,\mathcal{R},T}$ is the square root of the variance of the estimator $Q_{j,\mathcal{R};T}$ of $S_j(z)$, $z \in \mathcal{R}$. As in Spokoiny (1998), (4.4) can be viewed as a balance relation between the bias and the variance of this estimate. The $k_T$ term then appears as the correction term necessary in the pointwise estimation in order to bound the normalized loss [see Lepski (1990), Lepski and Spokoiny (1997)]. In the following results, we set $k_T$ proportional to $\log_2 T$.

PROPOSITION 5. *Let $\mathcal{R}$ be an interval of $(0,1)$ and consider the test rule (4.2). If the wavelet spectrum $S_j$ is regular on $\mathcal{R}$ in the sense of conditions (4.4)–(4.5), then, with $2\lambda = 2\eta = 2^{-j/2} 5(2\alpha + p)$ and $k_T \sim \log_2 T$,*

$$\Pr(\mathcal{R} \text{ is rejected}) = O(T^{-Kp\sqrt{\delta}})$$

*for some positive constant $K$ depending only on $K_2, K_3$ and $\|c\|_{1,\infty}$.*

We can also evaluate an upper bound for the $\ell_p$-risk associated with our estimator.

THEOREM 2. *Assume that the wavelet spectrum at scale $j$, $S_j(z)$, is homogeneous on the segment $\mathcal{R}$ in the sense of (4.4)–(4.5) with*

$$k_T \sim \log_2 T.$$

*If $\tilde{S}_j(z)$ is the pointwise estimator of the wavelet spectrum obtained by the estimation procedure (4.2)–(4.3) with*

$$\eta = 2^{-j/2} 5(2\alpha + p),$$

*then there exists $T_0$ such that the pointwise $\ell_p$-loss is bounded as follows:*

$$\mathrm{E} |\tilde{S}_j(z) - S_j(z)|^p \leq K \delta^{-p} T^{-p/2} (2^{j/2} \delta^{-1} + k_T)^p$$

*for $p \geq 2$ with a positive constant $K$ and $T > T_0$.*

The proof is found in Appendix B.8.

4.4. *Properties of the estimator in inhomogeneous regions.* We now describe the behavior of our estimator near a breakpoint located at a time point $z_\star$.

For a fixed scale $j \in \{-1, \ldots, -J_T\}$, assume the evolutionary wavelet spectrum to be homogeneous on $\mathcal{R}_0 = [z_0, z_\star)$ and on $\mathcal{R}_1 = (z_\star, z_1]$. We write $\mathcal{R} = \mathcal{R}_0 \cup \mathcal{R}_1 = [z_0, z_1]$ and

$$\theta_T := \mathrm{E}(Q_{j,\mathcal{R}_0;T} - Q_{j,\mathcal{R}_1;T})$$

and we assume that $\theta_T > 0$. The value of $\theta_T > 0$ precisely quantifies a change in the spectrum between the regions $\mathcal{R}_0$ and $\mathcal{R}_1$.

To prove the next proposition, we assume that the estimation procedure is such that $\mathcal{R}_0$ and $\mathcal{R}_1$ are in $\wp(\mathcal{R})$.

PROPOSITION 6. *If the evolutionary wavelet spectrum at scale $j$ contains a breakpoint at $z_\star$ (i.e., $\theta_T > 0$) and if $k_T \sim \log_2 T$, then*

$\Pr(\mathcal{R}$ is not rejected$)$

$$= O\bigg(\exp\bigg\{-\frac{T\theta_T^2(|\mathcal{R}_0|^2 \wedge |\mathcal{R}_1|^2)}{\log_2^2 T}\bigg\} + \exp\bigg\{-\frac{\sqrt{T}|\theta_T|(|\mathcal{R}_0| \wedge |\mathcal{R}_1|)}{\log_2^2 T}\bigg\}\bigg),$$

*where $c$ is a positive constant and $x \wedge y = \min(x, y)$.*

The proof of this proposition is given in Appendix B.9. Proposition 6 concerns the consistency of the test of homogeneity. Moreover, it allows a discussion of the local alternative to this test. We first note that the alternative hypothesis, that is, the definition of the inhomogeneous region, depends on the level of the jump $\theta_T$ and the lengths of the two segments $\mathcal{R}_0$ and $\mathcal{R}_1$. As a consequence, in order to study the local alternative, we need to investigate both cases $\theta_T \to 0$ and $(|\mathcal{R}_0| \wedge |\mathcal{R}_1|) \to 0$. It is interesting to note that Proposition 6 depends on the product $|\theta_T|(|\mathcal{R}_0| \wedge |\mathcal{R}_1|)$, and the local alternative of the test is then studied when this product tends to 0 as $T \to \infty$. From the proof of Proposition 6, it is straightforward to see that if

$$\frac{\log_2^2 T}{|\theta_T|(|\mathcal{R}_0| \wedge |\mathcal{R}_1|)\sqrt{T}} \to 0$$

as $T \to \infty$, then the estimation procedure is consistent in the sense that $\Pr(\mathcal{R}$ is not rejected$)$ is asymptotically zero.

**5. Simulation.** We conclude with a brief simulation study. We consider the evolutionary wavelet spectrum plotted in Figure 1 (upper plot). The first scale of this spectrum is given by $S_{-1}(z) = 1_{[0.25,0.575]}(z) + (\sin^2(2\pi z - \pi/4) + 0.5)1_{[0.75,1]}(z)$. The second scale is inactive. The other active scales are $S_{-3}(z) = (\sin(\pi z - \pi/4)^2 + 0.5)1_{[0,0.25]}(z)$ and $S_{-4}(z) = (\sin^2(5\pi z - \pi/4) +$ $0.5)1_{[0.375,1]}(z)$. We apply the estimation procedure to 100 different time series of length 1000 generated from this spectrum with Gaussian increments and Haar wavelets. For the sake of brevity, we only consider the estimation at the scale $j = -1$. The results of the 100 simulations are summarized in the upper plot of Figure 2. At each point of the 39 points of estimation, the vertical segment represents the median and the 90% interquantile interval from the 100 estimators. The bottom figure shows the estimator (bullet) from the single simulation given in Figure 1. The continuous line gives the estimator obtained from the ewspec function of the *WaveThresh 3* software package [Nason (1998)] using the recommendations suggested in this package for the choice of the parameters (other configurations performed quite similarly or worse). This estimator is a smoothing of the corrected wavelet periodogram using TI-wavelet soft thresholding; see Nason et al. (2000) for details. Note that this method is limited to dyadic sample sizes. As our simulation contains 1000 data, we repeat the last observation 24 times.

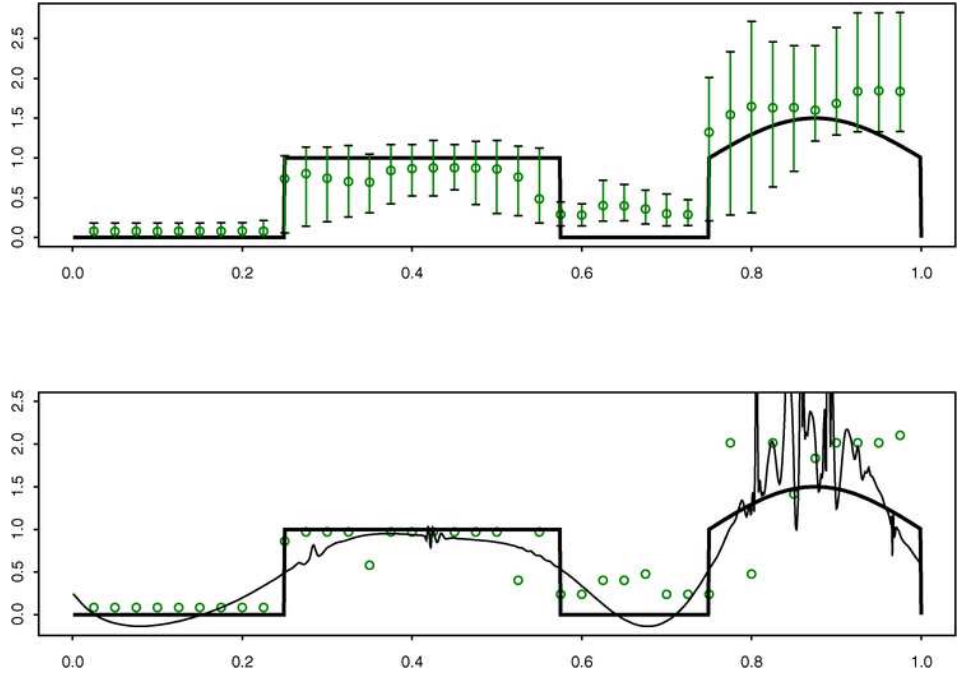

FIG. 2.  *The bold line in both graphs is the first scale of the evolutionary wavelet spectrum considered in Figure 1. The upper figure summarizes the results given from 100 simulations of the LSW process. In this figure, each vertical interval represents the 90% interquantile range from the 100 results and the bullet is the median. The bottom figure presents the local adaptive estimator (bullets) from the realization of the process shown in Figure 1 (lower plot). The continuous line is the estimator of Nason et al. (2000).*

The mean square error for the local adaptive estimator is lower (0.063) than for the nonlinear wavelet estimator (0.074). The mean absolute deviation is also lower (0.152 compared with 0.189 for the wavelet estimator). The lower plot of Figure 2 clearly shows the high variability of the `ewspec` estimator in the last part of the spectrum. We explain this phenomenon by the cross-correlation between the corrected wavelet periodograms at scales $-1$ and $-4$. It is interesting to note that our method seems to be more stable with respect to this phenomenon. This has been observed in comparison with `ewspec` using different wavelet families for smoothing.

In our simulation, it is worth mentioning that the local adaptive estimator is computed using the estimated variance, as explained in Section 3.3. Of course, there is a set of global parameters which must be chosen. For the example treated in this section, we set $M_T = 2$ and $|\mathcal{R}_T| = 9$ [see (3.8)]. With this, we have followed the guidelines given in the companion paper, Van Bellegem and von Sachs (2004) (Sections 2.3 and 2.4 therein) on the choice of nuisance parameters for the quadratic part of the estimator. In particular, two remaining global parameters have been chosen to equal the numerical values given for the (different) example of Section 2.5 therein. The paper also derives a new test of covariance stationarity and presents some applications to medical data analysis.

## APPENDIX A: PROPERTIES OF THE AUTOCORRELATION WAVELET SYSTEM

This section summarizes useful results on the system $\{\Psi_j\}$ and the operator $A$. Recall that we have denoted by $\mathcal{L}_j$ the length of $|\operatorname{supp}\psi_{j0}|$ for all $j = -1, -2, \ldots$, so we have $\mathcal{L}_j = (2^{-j} - 1)(\mathcal{L}_{-1} - 1) + 1 \leq 2^{-j}\mathcal{L}_{-1}$. We also recall the definition of the autocorrelation wavelet system $\{\Psi_j; j = -1, -2, \ldots\}$ which is the convolution of the nondecimated wavelet system,

$$\Psi_j(\tau) = \sum_{k=-\infty}^{\infty} \psi_{jk}(0)\psi_{jk}(\tau).$$

It is straightforward to check that $\Psi_j$ is compactly supported for all $j < 0$ and that the length of its support is bounded by $2\mathcal{L}_j - 1$.

The following lemma recalls other useful results on the autocorrelation wavelet system.

LEMMA A.1. (a) *For all scales $j$ and all $\tau$, $\Psi_j(\tau) = \Psi_j(-\tau)$.*

(b) *The autocorrelation wavelet system $\{\Psi_j; j = -1, -2, \ldots\}$ is linearly independent.*

(c) *The identity*

$$\sum_{j=-\infty}^{-1} 2^j \Psi_j(\tau) = \delta_0(\tau) \tag{A.1}$$

*holds for all $\tau \in \mathbb{Z}$.*

Property (a) is obvious and implies the symmetry of the local autocovariance function, that is, $c(z, \tau) = c(z, -\tau)$, as expected. Property (b) is proved as Theorem 1 of Nason et al. (2000) and shows that the local autocovariance function is univoquely defined. Finally, property (c) is proved as Lemma 6 of Fryźlewicz et al. (2003) and implies, for instance, that the wavelet spectrum of a white noise process is proportional to $2^j$ for all scales $j < 0$.

As the autocorrelation wavelet system is not orthogonal, we introduce the Gram matrix $A$ defined by $A_{j\ell} = \sum_\tau \Psi_j(\tau)\Psi_\ell(\tau)$. The following properties of $A$ are used thereafter.

LEMMA A.2. *For Haar and Shannon wavelets, there exists a finite positive constant $\nu$ such that the matrix $A$ fulfills the following properties for all $j = -1, \ldots, -\log_2 T$:*

$$\sum_{\ell=-\log_2 T}^{-1} A_{j\ell}^{-1} = 2^j + O(2^{j/2} T^{-1/2}); \tag{A.2}$$

$$\sum_{\ell=-\log_2 T}^{-1} |A_{j\ell}^{-1}| \leq \nu(1+\sqrt{2}) 2^{j/2}; \tag{A.3}$$

$$\sum_{\ell=-\log_2 T}^{-1} 2^{-\ell/2} |A_{j\ell}^{-1}| \leq \nu \cdot 2^{j/2} \log_2 T;$$

$$\sum_{\ell=-\log_2 T}^{-1} 2^{-\ell} |A_{j\ell}^{-1}| \leq \nu(2+\sqrt{2}) 2^{j/2} T^{1/2}. \tag{A.4}$$

*For all compactly supported wavelets, the matrix $A$ fulfills the following property:*

$$A_{j\ell} \leq (2\mathcal{L}_j - 1) \wedge (2\mathcal{L}_\ell - 1) \wedge \sqrt{\mathcal{L}_\ell \mathcal{L}_m}, \tag{A.5}$$

*where $x \wedge y = \min(x, y)$.*

PROOF. The following argument shows that the main term in (A.2) is $2^j$. Using the fact that $\Psi_\ell(0) = 1$ for all $\ell < 0$ and the identity (A.1), we may

write

$$\sum_{\ell=-\infty}^{-1} A_{j\ell}^{-1} = \sum_{\ell=-\infty}^{-1} A_{j\ell}^{-1} \sum_{m,u=-\infty}^{\infty} 2^m \Psi_m(u)\Psi_\ell(u)$$

$$= \sum_{m=-\infty}^{-1} 2^m \delta_0(j-m) = 2^j$$

from the definition of $A$. Observe that this argument holds for all compactly supported wavelets. To compute the remainder of (A.2), we introduce the auxiliary matrix $\Gamma = D' \cdot A \cdot D$ with diagonal matrix $D = \text{diag}(2^{\ell/2})_{\ell<0}$, that is, $\Gamma_{j\ell} = 2^{j/2} A_{j\ell} 2^{\ell/2}$. Nason et al. (2000), Theorem 2, have proven that the spectral norm of $\Gamma^{-1}$ is bounded for Haar and Shannon wavelets. We then get

$$\sum_{\ell=-\infty}^{-\log_2(T)-1} A_{j\ell}^{-1} = 2^{j/2} \sum_{\ell=-\infty}^{-\log_2(T)-1} 2^{\ell/2} \Gamma_{j\ell}^{-1} = O(2^{j/2} T^{-1/2}).$$

To prove (A.3), note that $\sum_{\ell=-\log_2 T}^{-1} |A_{j\ell}^{-1}| = \sum_{\ell=-\log_2 T}^{-1} 2^{j/2} 2^{\ell/2} |\Gamma_{j\ell}^{-1}| \leq 2^{j/2} \times (1+\sqrt{2})\nu$, using $\sup_{j\ell} |\Gamma_{j\ell}^{-1}| \leq \nu$. (A.4) is obtained similarly, using the approximation $\sum_{j=-\log_2 T}^{-1} 2^{-j/2} \leq (2+\sqrt{2})\sqrt{T}$. (A.5) follows from the definition of $A_{j\ell}$ and the support of the autocorrelation wavelets, using $|\Psi_j(\tau)| \leq 1$ uniformly in $j$ and $\tau$. $\square$

## APPENDIX B: PROOFS

Suppose $M$ is an $n \times n$ matrix and $\overline{M}'$ is the conjugate transpose of $M$. We denote by

$$\|M\|_2 := \sqrt{\text{tr}(\overline{M}'M)}$$

the Euclidean norm of $M$ and by

$$\|M\|_{\text{spec}} := \max\{\sqrt{\lambda} : \lambda \text{ is an eigenvalue of } M^\star M\}$$

the spectral norm of $M$. If $M$ is symmetric and nonnegative definite, then by standard theory, we have $\|M\|_{\text{spec}} = \sup\{\|Mx\|_2 : x \in \mathbb{C}^n, \|x\|_2 = 1\}$. We will also use the following standard relations which hold for all symmetric matrices $B, C$:

(B.1) $\quad \|B\|_{\text{spec}} \leq \|B\|_2;$

(B.2) $\quad \|B\|_{\text{spec}} = \max\{\lambda : \lambda \text{ is an eigenvalue of } B\};$

(B.3) $\quad \|BC\|_{\text{spec}} \leq \|B\|_{\text{spec}} \|C\|_{\text{spec}};$

(B.4) $\quad \|BC\|_2 \leq \|B\|_{\text{spec}} \|C\|_2 \leq \|B\|_2 \|C\|_2.$

In the sequel, we use the convention $w_{jk;T} = 0$ for $k < 0$ and $k \geq T$, which leads to helpful simplifications in the following proofs.

**B.1. Proof of Proposition 1.** On one hand, due to Definition 1 and equation (2.2), we have

$$c_{X,T}(z,\tau) = \text{Cov}(X_{[zT],T}, X_{[zT]+\tau,T}) = \sum_{j=-\infty}^{-1} \sum_{k=-\infty}^{\infty} |w_{j,k+[zT];T}|^2 \psi_{jk}(0)\psi_{jk}(\tau)$$

$$= \sum_{j=-\infty}^{-1} \sum_{k=-\infty}^{\infty} S_j\left(\frac{k+[zT]}{T}\right)\psi_{jk}(0)\psi_{jk}(\tau) + \text{Rest}_T(z,\tau),$$

where the remainder is such that

$$|\text{Rest}_T(z,\tau)| = O(T^{-1}) \sum_{j=-\infty}^{-1} \sum_{k=-\infty}^{\infty} C_j |\psi_{jk}(0)\psi_{jk}(\tau)|,$$

by assumption (2.2). On the other hand, we have $c_X(z,\tau) = \sum_{j=-\infty}^{-1} \sum_{k=-\infty}^{\infty} S_j(z)\psi_{jk}(0)\psi_{jk}(\tau)$. Then,

$$\sum_{\tau=-\infty}^{\infty} \int_0^1 dz |c_{X,T}(z,\tau) - c_X(z,\tau)|$$

$$\leq \sum_{\tau=-\infty}^{\infty} \int_0^1 dz \sum_{j=-\infty}^{-1} \sum_{k=-\infty}^{\infty} \left|S_j\left(\frac{k+[zT]}{T}\right) - S_j(z)\right| |\psi_{jk}(0)\psi_{jk}(\tau)|$$

$$+ \sum_{\tau=-\infty}^{\infty} \int_0^1 dz |\text{Rest}_T(z,\tau)|.$$

With appropriate changes of variable, this bound may be written as

$$\sum_{\tau=-\infty}^{\infty} \sum_{j=-\infty}^{-1} \sum_{k=-\infty}^{\infty} \sum_{t=0}^{T-1} \int_0^{1/T} dz \left|S_j\left(\frac{k+[zT]+t}{T}\right) - S_j\left(z+\frac{t}{T}\right)\right| |\psi_{jk}(0)\psi_{jk}(\tau)|$$

$$+ \sum_{\tau=-\infty}^{\infty} \int_0^1 dz |\text{Rest}_T(z,\tau)|,$$

which is bounded by

$$T^{-1} \sum_{\tau=-\infty}^{\infty} \sum_{j=-\infty}^{-1} \sum_{k=-\infty}^{\infty} |k| \text{TV}(S_j) |\psi_{jk}(0)\psi_{jk}(\tau)| + \sum_{\tau=-\infty}^{\infty} \int_0^1 dz |\text{Rest}_T(z,\tau)|,$$

where we have used the following property of the total variation:

$$\sum_{t=0}^{T-1} \left|S_j\left(\frac{t}{T}+\frac{\alpha}{T}\right) - S_j\left(\frac{t}{T}+\frac{\beta}{T}\right)\right| \leq |\alpha-\beta| \text{TV}(S_j) \qquad \text{for all } \alpha,\beta \in \mathbb{N}.$$

(B.5)

As the support of $\psi_{jk}(0)$ is of length $\mathcal{L}_j$, we get $|k| \leq \mathcal{L}_j$ in the first term. Together with condition (2.3) of Definition 1, this finally leads to

$$\sum_{\tau=-\infty}^{\infty} \int_0^1 dz |c_{X,T}(z,\tau) - c_X(z,\tau)|$$

$$\leq O(T^{-1}) \sum_{j=-\infty}^{-1} (C_j + \mathcal{L}_j L_j) \sum_{\tau,k=-\infty}^{\infty} |\psi_{jk}(0)\psi_{jk}(\tau)|.$$

The compact support of $\psi_{jk}$ limits the sums over $k$ and $\tau$ as follows:

$$(B.6) \quad \sum_{\tau,k=-\infty}^{\infty} |\psi_{jk}(0)\psi_{jk}(\tau)| = \sum_{\tau=-\mathcal{L}_j+1}^{\mathcal{L}_j-1} \sum_{k=-\infty}^{\infty} |\psi_{jk}(0)\psi_{jk}(\tau)| \leq 2\mathcal{L}_j - 1,$$

by the Cauchy–Schwarz inequality for the sum over $k$. We then get the result by (2.4).

**B.2. Preliminary results.** Define $\underline{X}_T = (X_{0,T}, \ldots, X_{T-1,T})'$. By definition, $Q_{j,\mathcal{R};T}$ can be decomposed into the sum of a quadratic and a linear form,

$$(B.7) \qquad Q_{j,\mathcal{R};T} = Q^{\circ}_{j,\mathcal{R};T} + q^{\circ}_{j,\mathcal{R};T},$$

where

$$(B.8) \qquad Q^{\circ}_{j,\mathcal{R};T} = \underline{X}'_T U_{j,\mathcal{R};T} \underline{X}_T$$

is a quadratic form with the $T \times T$ matrix $U_{j,\mathcal{R};T}$ whose entry $(s,t)$ is

$$U_{st} = |\mathcal{R}T|^{-1} \sum_{\ell=-\log_2 T}^{-1} A_{j\ell}^{-1} \sum_{k \in \mathcal{R}T} \psi_{\ell k}(s) \psi_{\ell k}(t)$$

and $q^{\circ}_{j,\mathcal{R};T} = |\mathcal{R}T|^{-1} \sum_{k \in \mathcal{R}T} z_{j,k;T}$ is the linear form. For notational convenience, we omit the dependence of $U_{st}$ on $j$ and $\mathcal{R}$. Assuming that the orthonormal increment processes $\{\xi_{jk}\}$ in Definition 1 are Gaussian, $\underline{X}_T$ is a multivariate Gaussian random variable with covariance matrix $\Sigma_T = \mathrm{Cov}(\underline{X}_T \underline{X}'_T)$. Therefore, we can write

$$Q_{j,\mathcal{R};T} = \underline{Z}'_T M_{j,\mathcal{R};T} \underline{Z}_T + q^{\circ}_{j,\mathcal{R};T},$$

where $\underline{Z}_T = (Z_1, \ldots, Z_T)'$ is a vector of i.i.d. Gaussian random variables with zero mean and $\mathrm{Var}\, Z_1 = 1$, and

$$(B.9) \qquad M_{j,\mathcal{R};T} = \Sigma_T'^{1/2} U_{j,\mathcal{R};T} \Sigma_T^{1/2}$$

is the matrix of the quadratic form.

In our proofs, we use the following lemma quoted from Neumann (1996).

LEMMA B.1. *Let $\underline{Z}_n = (Z_1, \ldots, Z_n)'$ be a vector of i.i.d. Gaussian random variables with zero mean and $\mathrm{Var}\, Z_1 = 1$. If $M_n$ is an $n \times n$ real matrix, then*

$$\mathrm{E}(\underline{Z}_n' M_n \underline{Z}_n) = \mathrm{tr}\, M_n,$$

$$\mathrm{Var}(\underline{Z}_n' M_n \underline{Z}_n) = 2 \mathrm{tr}\, M_n' M_n = 2\|M_n\|_2^2$$

*and, for all $r \geq 2$, if $\mathrm{Cum}_r$ denotes the $r$th cumulant, we have*

$$|\mathrm{Cum}_r(\underline{Z}_n' M_n \underline{Z}_n)| \leq 2^{r-1}(r-1)! \|M_n\|_2^2 \{\lambda_{\max}(M_n)\}^{r-2}.$$

The following lemmas derive some bounds for the Euclidean and spectral norms of $U_{j,\mathcal{R};T}$ and $\Sigma_T$.

LEMMA B.2. *With fixed $\mathcal{R} \subseteq (0,1)$, there exists a $T_0$ such that, uniformly in $T \geq T_0$,*

$$\|U_{j,\mathcal{R};T}\|_2^2 \leq K_2^2 2^j |\mathcal{R}|^{-2} T^{-1}$$

*for all $j = -1, \ldots, J_T = o_T(\log_2 T)$, where $K_2$ depends only on the mother wavelet $\psi$.*

PROOF. If we let $\mathcal{R} = (r_1, r_2) \subseteq (0,1)$, then we can write $U_{st} = U_{st}^{(2)} - U_{st}^{(1)}$, where $U_{st}^{(1)} := |\mathcal{R}T|^{-1} \sum_\ell A_{j\ell}^{-1} \sum_{k=0}^{[r_1 T]-1} \psi_{\ell k}(t) \psi_{\ell k}(s)$ is the element $(s,t)$ of a matrix $U_{j,\mathcal{R};T}^{(1)}$ and $U_{st}^{(2)} := |\mathcal{R}T|^{-1} \sum_\ell A_{j\ell}^{-1} \sum_{k=0}^{[r_2 T]} \psi_{\ell k}(t) \psi_{\ell k}(s)$ is the element $(s,t)$ of a matrix $U_{j,\mathcal{R};T}^{(2)}$. Note that the compact support of the wavelet $\psi$ implies that $U_{st}^{(1)} = 0$ when $s$ or $t \geq [r_1 T]$ and, similarly, $U_{st}^{(2)} = 0$ when $s$ or $t > [r_2 T]$. We also introduce the matrix $U_{j,\mathcal{R};T}^{\star(1)}$ whose entry $(s,t)$ is $U_{st}^{\star(1)} := |\mathcal{R}T|^{-1} \sum_\ell A_{j\ell}^{-1} \Psi_\ell(s-t) \mathbb{I}_{0 \leq s,t < [r_1 T]}$ and similarly define $U_{j,\mathcal{R};T}^{\star(2)}$. We now have the decomposition

$$\|U_{j,\mathcal{R};T}\|_2^2 \leq 2\|U_{j,\mathcal{R};T}^{(1)} - U_{j,\mathcal{R};T}^{(1)\star}\|_2^2 + 4\|U_{j,\mathcal{R};T}^{(2)} - U_{j,\mathcal{R};T}^{(2)\star}\|_2^2$$
$$+ 4\|U_{j,\mathcal{R};T}^{(1)\star} - U_{j,\mathcal{R};T}^{(2)\star}\|_2^2.$$

From the definition of the autocorrelation wavelet $\Psi$, the first term is

$$\|U_{j,\mathcal{R};T}^{(1)} - U_{j,\mathcal{R};T}^{(1)\star}\|_2^2$$
$$= |\mathcal{R}T|^{-2} \sum_{\ell,m=-\log_2 T}^{-1} A_{j\ell}^{-1} A_{jm}^{-1} \sum_{s,t=0}^{[r_1 T]-1} \sum_{k,n=[r_1 T]}^{\infty} \psi_{\ell k}(t) \psi_{\ell k}(s) \psi_{mn}(t) \psi_{mn}(s).$$

The compact support of $\psi_{\ell k}(s)$ implies that $s > k - \mathcal{L}_\ell \geq ([r_1 T] - \mathcal{L}_\ell) \vee 0$. Using the same argument on $\psi_{mn}(t)$, we have $t > ([r_1 T] - \mathcal{L}_m) \vee 0$. Using the Cauchy–Schwarz inequality twice for the sums over $k$ and $n$, we get the bound

$$\|U_{j,\mathcal{R};T}^{(1)} - U_{j,\mathcal{R};T}^{(1)\star}\|_2^2 \leq |\mathcal{R}T|^{-2}\left(\sum_{\ell=-\log_2 T}^{-1} \mathcal{L}_\ell |A_{j\ell}^{-1}|\right)^2$$

$$\leq |\mathcal{R}T|^{-2}\nu^2(2+\sqrt{2})^2 2^j T \mathcal{L}_{-1}^2,$$

by (A.4). The second term is bounded similarly. The third term is bounded by $2\|U_{j,\mathcal{R};T}^{(1)\star}\|_2^2 + 2\|U_{j,\mathcal{R};T}^{(2)\star}\|_2^2$ and each term of this last sum can be bounded as

$$\|U_{j,\mathcal{R};T}^{(1)\star}\|_2^2$$

$$\leq |\mathcal{R}T|^{-2} \sum_{s=0}^{T-1} \sum_{t=-\infty}^{\infty} \sum_{\ell,m} A_{k\ell}^{-1} A_{jm}^{-1} \Psi_\ell(s-t)\Psi_m(s-t)$$

$$= T|\mathcal{R}T|^{-2} A_{jj}^{-1},$$

which leads to the result. $\square$

Finally, the proof of the following lemma is similar to the proof of Lemma 5.9 in Dahlhaus and Polonik (2006).

LEMMA B.3. *Under assumption (3.5), $\|\Sigma_T\|_{\mathrm{spec}} = \|\Sigma_T^{1/2}\|_{\mathrm{spec}}^2 \leq \|c_X\|_{1,\infty} < \infty$.*

**B.3. Proof of Proposition 2.**

*Expectation.* In decomposition (B.7), we first note that $\mathrm{E}\,q_{j,\mathcal{R};T}^\circ = 0$. Next, a straightforward expansion leads to

$$\mathrm{E}\,Q_{j,\mathcal{R};T}^\circ = |\mathcal{R}T|^{-1} \sum_{k\in\mathcal{R}T} \sum_{\ell=-\log_2 T}^{-1} A_{j\ell}^{-1} \sum_{s,t=0}^{T-1} \psi_{\ell k}(s)\psi_{\ell k}(t)$$

$$\times \sum_{m=-\infty}^{-1} \sum_{n=-\infty}^{\infty} w_{mn;T}^2 \psi_{mn}(s)\psi_{mn}(t)$$

$$= |\mathcal{R}T|^{-1} \sum_{k\in\mathcal{R}T} \sum_{\ell=-\log_2 T}^{-1} A_{j\ell}^{-1}$$

$$\times \sum_{m=-\infty}^{-1} \sum_{n=-\infty}^{\infty} w_{mn;T}^2 \left(\sum_{s=0}^{T-1} \psi_{\ell k}(s)\psi_{mn}(s)\right)^2.$$

Defining $u := n - k$, we can write

$$\mathrm{E}\, Q^{\circ}_{j,\mathcal{R};T} = |\mathcal{R}T|^{-1} \sum_{k \in \mathcal{R}T} \sum_{m=-\infty}^{-1} \sum_{u=-\infty}^{\infty} w^2_{m,u+k,T}$$

$$\times \sum_{\ell=-\log_2 T}^{-1} A^{-1}_{j\ell} \left( \sum_{s=-\infty}^{\infty} \psi_{\ell k}(s) \psi_{m,u+k}(s) \right)^2.$$

By Definition 1, we can write $w^2_{m,u+k,T} = S_m(k/T) + R_T(m,u,k)$ with

$$|R_T(m,u,k)| \leq \left| S_m\!\left(\frac{u+k}{T}\right) - S_m\!\left(\frac{k}{T}\right) \right| + \frac{\overline{C} C_m}{T},$$

which leads to

$$\mathrm{E}\, Q^{\circ}_{j,\mathcal{R};T} = |\mathcal{R}T|^{-1} \sum_{k \in \mathcal{R}T} \sum_{m=-\infty}^{-1} S_m\!\left(\frac{k}{T}\right)$$

$$\times \sum_{\ell=-\log_2 T}^{-1} A^{-1}_{j\ell} \sum_{u=-\infty}^{\infty} \left( \sum_{s=-\infty}^{\infty} \psi_{\ell k}(s) \psi_{m,u+k}(s) \right)^2 + \mathrm{Rest}_T.$$

By construction of the matrix $A$, we observe that

$$(\mathrm{B}.10) \qquad A_{\ell m} = \sum_{u=-\infty}^{\infty} \left( \sum_{s=-\infty}^{\infty} \psi_{\ell k}(s) \psi_{m,u+k}(s) \right)^2,$$

which implies, by Assumption 4, that

$$\mathrm{E}\, Q^{\circ}_{j,\mathcal{R};T} = |\mathcal{R}T|^{-1} \sum_{k \in \mathcal{R}T} S_j\!\left(\frac{k}{T}\right) + \mathrm{Rest}_T$$

(B.11)

$$= |\mathcal{R}|^{-1} \int_{\mathcal{R}} dz\, S_j(z) + O(|\mathcal{R}T|^{-1} L_j) + \mathrm{Rest}_T,$$

where the last equality is a standard result on the total variation [see, e.g., Lemma P5.1 of Brillinger (1975)].

We now bound $|\mathrm{Rest}_T|$. As $s$ goes from $-\infty$ to $\infty$, we have

$$|\mathrm{Rest}_T| \leq \sum_{m=-\infty}^{-1} \sum_{\ell=-\log_2 T}^{-1} |A^{-1}_{j\ell}|$$

$$\times \sum_{u=-\infty}^{\infty} |\mathcal{R}T|^{-1} \sum_{k \in \mathcal{R}T} \left\{ \left| S_m\!\left(\frac{u+k}{T}\right) - S_m\!\left(\frac{k}{T}\right) \right| + \frac{\overline{C} C_m}{T} \right\}$$

$$\times \left( \sum_{s=-\infty}^{\infty} \psi_{\ell 0}(s) \psi_{mu}(s) \right)^2.$$

Using (B.5) for the sum over $k$, $|\mathrm{Rest}_T|$ is bounded by

$$\sum_{m=-\infty}^{-1} \sum_{u=-\infty}^{\infty} \left\{|u|\frac{\mathrm{TV}(S_m)}{|\mathcal{R}T|} + \frac{\overline{C}C_m}{T}\right\} \sum_{\ell=-\log_2 T}^{-1} |A_{j\ell}^{-1}| \left(\sum_{s=-\infty}^{\infty} \psi_{\ell 0}(s)\psi_{mu}(s)\right)^2.$$

In this last expression, the compact support of $\psi_{\ell 0}$ and $\psi_{mu}$ implies that $|u| \leq \mathcal{L}_\ell \vee \mathcal{L}_m$, where $x \vee y = \max(x,y)$. Together with (B.10), we get

$$|\mathrm{Rest}_T| \leq |\mathcal{R}T|^{-1} \sum_{m=-\infty}^{-1} \sum_{\ell=-\log_2 T}^{-1} \{\mathrm{TV}(S_m)(\mathcal{L}_\ell \vee \mathcal{L}_m) + \overline{C}C_m\}|A_{j\ell}^{-1}|A_{\ell m},$$

which, with (A.5), leads to

$$\begin{aligned}
|\mathrm{Rest}_T| \\
\leq |\mathcal{R}T|^{-1} \sum_{m,\ell} &\{\mathrm{TV}(S_m)\mathcal{L}_\ell(2\mathcal{L}_m - 1) \\
&+ \mathrm{TV}(S_m)\mathcal{L}_m(2\mathcal{L}_\ell - 1) + \overline{C}C_m(2\mathcal{L}_m - 1)\}|A_{j\ell}^{-1}| \\
= 2(2+\sqrt{2})&\nu 2^{j/2}|\mathcal{R}T|^{-1}\sqrt{T}\mathcal{L}_{-1} \\
\times \sum_{m=-\infty}^{-1}& (2\mathcal{L}_m - 1)\mathrm{TV}(S_m) + O(2^{j/2}|\mathcal{R}T|^{-1}),
\end{aligned}$$
(B.12)

using (A.4) and (2.4).

*Variance.* Using decomposition (B.7), the variance is decomposed as $\mathrm{Var}\, Q_{j,\mathcal{R};T} = \mathrm{Var}\, Q_{j,\mathcal{R};T}^\circ + \mathrm{Var}\, q_{j,\mathcal{R};T}^\circ$, where $\mathrm{Var}\, q_{j,\mathcal{R};T}^\circ = C^2 2^j/|\mathcal{R}T|$. Using Lemma B.1 with (B.4), we can write $\mathrm{Var}\, Q_{j,\mathcal{R};T}^\circ = 2\|M_{j,\mathcal{R};T}\|_2^2 \leq 2\|\Sigma_T^{1/2}\|_{\mathrm{spec}}^4 \times \|U_{j,\mathcal{R};T}\|_2^2$ and the result follows from Lemmas B.2 and B.3.

**B.4. Proof of Proposition 3 and its consequences.** Our proof of Proposition 3 requires the use of an exponential bound for linear and quadratic forms of Gaussian random variables. For the sake of presentation, we here summarize the results we use.

PROPOSITION B.1. *Let $Z$ be a Gaussian random variable with mean zero and unit variance. Then, for all $\lambda > 0$,*

$$\Pr(|Z| > \lambda) \leq \left(1 \wedge \frac{1}{\lambda\sqrt{2\pi}}\right)e^{-\lambda^2/2},$$

*where $a \wedge b = \min(a,b)$.*

Let $\underline{Z}_n = (Z_1, \ldots, Z_n)'$ be a vector of i.i.d. Gaussian random variables with zero mean and $\operatorname{Var} Z_1 = 1$. If $M_n$ is an $n \times n$ matrix such that $\|M_n\|_{\operatorname{spec}} \leq \tau_\infty$ and $\sigma_n^2 = 2\|M_n\|_2^2$, then, for all $\lambda > 0$,

$$\Pr((\underline{Z}_n' M_n \underline{Z}_n - \operatorname{tr} M_n) > \sigma_n \lambda) \leq 2 \exp\left(-\frac{1}{4} \cdot \frac{\lambda^2}{1 + 2(\lambda \tau_\infty / \sigma_n)}\right).$$

Moreover, if $Y$ is a Gaussian random variable with mean zero and variance $\sigma^2 \leq \sigma_n^2$, then

$$\Pr((\underline{Z}_n' M_n \underline{Z}_n + Y - \operatorname{tr} M_n) > \sigma_n \lambda) \leq 3 \exp\left(-\frac{1}{4} \cdot \frac{\lambda^2}{1 + 2(\lambda \tau_\infty / \sigma_n)}\right).$$

PROOF. We prove the first inequality. On the one hand, by Chebyshev's inequality,

$$\Pr(Z > \lambda) \leq \inf_{t>0} \exp\{-t\lambda + \log \operatorname{E}(e^{tZ})\},$$

where $\operatorname{E}(e^{tZ}) = e^{-t^2/2}$. The minimum is attained for $t = \lambda$ and we get $\Pr(|Z| > \lambda) \leq e^{-\lambda^2/2}$. On the other hand, a straightforward calculation leads to

$$\Pr(Z > \lambda) = \int_\lambda^\infty \frac{1}{\sqrt{2\pi}} e^{-t^2/2} \, dt \leq \int_\lambda^\infty \frac{\lambda}{\sqrt{2\pi}} e^{-t^2/2} \, dt = \frac{1}{\lambda\sqrt{2\pi}} e^{-\lambda^2/2}$$

and the result follows. The second inequality follows the proof of Proposition A.1 in Dahlhaus and Polonik (2006). The last inequality is derived from the two former inequalities. □

As in the proof of Proposition 2, equation (B.9), we write $Q_{j,\mathcal{R};T}$ as a quadratic form of Gaussian variables in order to apply Proposition B.1 with $M_{j,\mathcal{R};T} = \Sigma_T'^{1/2} U_{j,\mathcal{R};T} \Sigma_T^{1/2}$ and thereby prove the assertion.

PROOF OF PROPOSITION 3. We use the last exponential inequality of Proposition B.1 because $Q_{j,\mathcal{R};T}$ can be decomposed [see (B.7)] into $Q_{j,\mathcal{R};T}^\circ + q_{j,\mathcal{R};T}^\circ$, where $Q_{j,\mathcal{R};T}^\circ = \underline{Z}_T' M_{j,\mathcal{R};T} \underline{Z}_T$ and $q_{j,\mathcal{R};T}^\circ \sim \mathcal{N}(0, C^2 2^j / |\mathcal{R}T|)$. Note that Lemmas B.2 and B.3 imply, with (B.1) and (B.3), that

(B.13) $$\|M_{j,\mathcal{R};T}\|_{\operatorname{spec}} \leq 2^{j/2} K_2 \|c_X\|_{1,\infty} |\mathcal{R}|^{-1} T^{-1/2}.$$

Therefore, Proposition B.1 leads to

$$\Pr((Q_{j,\mathcal{R};T} - Q_{j,\mathcal{R}}) > \eta \sigma_{j,\mathcal{R},T})$$
$$\leq \Pr((Q_{j,\mathcal{R};T} - \operatorname{E} Q_{j,\mathcal{R};T}) > \eta \sigma_{j,\mathcal{R},T}/2)$$
$$+ \exp\left(1 - \frac{\eta \sigma_{j,\mathcal{R},T}}{2|\operatorname{E} Q_{j,\mathcal{R};T} - Q_{j,\mathcal{R}}|}\right)$$

$$\leq 3\exp\left(-\frac{1}{16}\cdot\frac{\eta^2}{1+\eta(2^{j/2}K_2\|c_X\|_{1,\infty})/(|\mathcal{R}|T^{1/2}\sigma_{j,\mathcal{R},T})}\right)$$
$$+\exp\left(1-\frac{\eta\sigma_{j,\mathcal{R},T}}{2|\operatorname{E}Q_{j,\mathcal{R};T}-Q_{j,\mathcal{R}}|}\right).$$

To bound the second probability, we observe that (B.11) and (B.12) lead to $|\operatorname{E}Q_{j,\mathcal{R};T}-Q_{j,\mathcal{R}}|\leq|\mathcal{R}T|^{-1}(L_j+K_3 2^{(j/2)-1}\sqrt{T})$ with $K_3=4\nu(2+\sqrt{2})(2\rho-1)(\overline{C}\vee 1)\mathcal{L}_{-1}$. This implies

$$\Pr((Q_{j,\mathcal{R};T}-Q_{j,\mathcal{R}})\geq\eta\sigma_{j,\mathcal{R},T})$$
$$\leq 3\exp\left(-\frac{1}{16}\cdot\frac{\eta^2\sigma_{j,\mathcal{R},T}}{\sigma_{j,\mathcal{R},T}+\eta(2^{j/2}K_2\|c_X\|_{1,\infty}\sqrt{T})/|\mathcal{R}T|}\right)$$
$$+\exp\left(1-\frac{1}{2\eta}\frac{\eta^2\sigma_{j,\mathcal{R},T}}{(L_j+K_3 2^{(j/2)-1}\sqrt{T})/|\mathcal{R}T|}\right)$$

and the result follows. $\square$

PROOF OF COROLLARY 1. In the following proof, $K$ denotes a generic constant and $k_T$ is an increasing function of $T$. By Proposition 2, $\sigma^2_{j,\mathcal{R},T}:=\operatorname{Var}Q_{j,\mathcal{R};T}\leq(C^2+c^2/|\mathcal{R}|)2^j/|\mathcal{R}T|$ uniformly in $j$, which implies

$$\Pr\left(\sup_{-J_T\leq j<0}|Q_{j,\mathcal{R};T}-Q_{j,\mathcal{R}}|\geq k_T\sqrt{(C^2+c^2/|\mathcal{R}|)/|\mathcal{R}T|}\right)$$
$$\leq\sum_{j=-J_T}^{-1}\Pr(|Q_{j,\mathcal{R};T}-Q_{j,\mathcal{R}}|\geq 2^{-j/2}k_T\sigma_{j,\mathcal{R},T}).$$

Using Proposition 3, this probability is bounded by

$$c_0 J_T\max_{-J_T\leq j<0}\exp\Bigg\{-\frac{1}{16}\cdot(2^{-j}k_T^2/2)\Big/\bigg[1+\frac{2k_T 2^{-j/2}L_j}{(|\mathcal{R}T|\sigma_{j,\mathcal{R},T})}$$
$$+\frac{k_T\sqrt{T}}{|\mathcal{R}T|\sigma_{j,\mathcal{R},T}}(K_2\|c_X\|_{1,\infty}+K_3)\bigg]\Bigg\}.$$

Proposition 2 shows that, for $T$ sufficiently large, $\sigma_{j,\mathcal{R},T}\geq\sqrt{2^j/|\mathcal{R}T|}$. This leads to the bound

$$c_0 J_T\max_{-J_T\leq j<0}\exp\Bigg\{-\frac{1}{16}\cdot(k_T^2/2)\Big/\bigg[2^j+\frac{k_T 2^{-j/2}L_j}{\sqrt{|\mathcal{R}T|}}$$
$$+\frac{2^{j/2}k_T}{\sqrt{|\mathcal{R}|}}(K_2\|c_X\|_{1,\infty}+K_3)\bigg]\Bigg\}.$$

By assumption (2.4), there exists a positive constant $\rho'$ such that $L_j \leq 2^{j/2}\rho'$. Then, asymptotically, the rate of convergence of the dominant terms in this exponential expression are given by $J_T \cdot \exp(-k_T)$, which is $o_T(1)$ by the assumption on $k_T$. $\square$

### B.5. Proof of Proposition 4.

LEMMA B.4. *If $U_{ts}^{(j)} = |\mathcal{R}T|^{-1} \sum_{\ell=-\log_2 T}^{-1} A_{j\ell}^{-1} \sum_{k\in\mathcal{R}T} \psi_{\ell k}(s)\psi_{\ell k}(t)$, then*

$$\sum_{t=-\infty}^{\infty} \sum_{s,u=-\infty}^{\infty} |U_{ts}^{(j)} U_{tu}^{(j)}| \mathbb{I}_{|s-u|\leq N_T} \leq 2^{j+1}\mathcal{L}_{-1}\nu^2 \frac{TN_T \log_2^2 T}{|\mathcal{R}T|^2}$$

$$= O\left(2^j \frac{N_T \log_2^2 T}{T}\right).$$

PROOF. Direct calculation yields

$$\sum_{t=-\infty}^{\infty} \sum_{s,u=-\infty}^{\infty} |U_{ts}^{(j)} U_{tu}^{(j)}| \mathbb{I}_{|s-u|\leq N_T}$$

$$\leq |\mathcal{R}T|^{-2} \sum_{\ell,m=-\log_2 T}^{-1} |A_{j\ell}^{-1}||A_{jm}^{-1}| \sum_{s,u=-\infty}^{\infty} \mathbb{I}_{|s-u|\leq N_T}$$

$$\times \sum_{t=-\infty}^{\infty} \left(\sum_{k\in\mathcal{R}T} |\psi_{\ell k}(s)\psi_{\ell k}(t)|\right)\left(\sum_{n\in\mathcal{R}T} |\psi_{mn}(u)\psi_{mn}(t)|\right).$$

Using the Cauchy–Schwarz inequality for the sum over $t$, we get a product of two terms similar to $(\sum_t(\sum_k \psi_{\ell k}(s)\psi_{\ell k}(t))^2)^{1/2} \leq \sqrt{2\mathcal{L}_\ell - 1}$. Then,

$$\sum_{t=-\infty}^{\infty} \sum_{s,u=-\infty}^{\infty} |U_{ts}^{(j)} U_{tu}^{(j)}| \mathbb{I}_{|s-u|\leq N_T}$$

$$\leq TN_T|\mathcal{R}T|^{-2} \sum_{\ell,m} |A_{j\ell}^{-1}||A_{jm}^{-1}|\sqrt{2\mathcal{L}_\ell - 1}\sqrt{2\mathcal{L}_m - 1}$$

and we obtain the result by (A.4). $\square$

In the proof of Proposition 4, we need a modification of Corollary 1, in which $\mathcal{R}$ is replaced by $\mathcal{R}_T$. The proof of the following result is along the lines of the proof of Corollary 1.

LEMMA B.5. *Under the assumptions of Propositions 2 and 3, there exists $T_0 \geq 1$ such that, for all $T \geq T_0$,*

$$\Pr\left(\sup_{-J_T \leq j < 0} |Q_{j,\mathcal{R}_T(s);T} - Q_{j,\mathcal{R}_T(s)}| \geq \frac{k_T}{|\mathcal{R}_T|}\sqrt{\frac{C^2+c^2}{T}}\right) = o_T(1),$$

*provided that* $J_T \cdot \exp(-k_T\sqrt{|\mathcal{R}_T|}) = o_T(1)$.

PROOF OF PROPOSITION 4. Define $\bar{\sigma}_{s,s+u} := \sum_{\ell=-\log_2 T}^{-1} Q_{\ell,\mathcal{R}_T(s)} \Psi_\ell(u) \times \mathbb{I}_{|u| \leq M_T}$, the entries of a matrix $\bar{\Sigma}$, and define $\bar{\sigma}_{j,\mathcal{R},T}^2 := 2\|U'_{j,\mathcal{R};T}\bar{\Sigma}_T\|_2^2 + C^2 2^j/|\mathcal{R}T|$. Our proof is based on the decomposition

$$\tilde{\sigma}_{j,\mathcal{R},T}^2 - \sigma_{j,\mathcal{R},T}^2 = (\tilde{\sigma}_{j,\mathcal{R},T}^2 - \bar{\sigma}_{j,\mathcal{R},T}^2) + (\bar{\sigma}_{j,\mathcal{R},T}^2 - \sigma_{j,\mathcal{R},T}^2),$$

where the first term is stochastic while the second term is deterministic.

We will first show that the deterministic term $|\bar{\sigma}_{j,\mathcal{R},T}^2 - \sigma_{j,\mathcal{R},T}^2|$ is $o(2^{j-J_T}T^{-1})$. Using (B.4), we can write

$$\frac{1}{2}(\bar{\sigma}_{j,\mathcal{R},T}^2 - \sigma_{j,\mathcal{R},T}^2)$$
$$= \|U'_{j,\mathcal{R};T}\bar{\Sigma}_T\|_2^2 - \|U'_{j,\mathcal{R};T}\Sigma_T\|_2^2$$
$$\leq \|U'_{j,\mathcal{R};T}(\bar{\Sigma}_T - \Sigma_T)\|_2^2 + 2 \cdot \|U'_{j,\mathcal{R};T}\Sigma_T\|_2 \cdot \|U'_{j,\mathcal{R};T}(\bar{\Sigma}_T - \Sigma_T)\|_2$$
$$\leq \|U_{j,\mathcal{R};T}\|_2^2 \cdot \|\bar{\Sigma}_T - \Sigma_T\|_{\mathrm{spec}}^2$$
$$\quad + 2 \cdot \|U_{j,\mathcal{R};T}\|_2^2 \cdot \|\Sigma_T\|_{\mathrm{spec}} \cdot \|\bar{\Sigma}_T - \Sigma_T\|_{\mathrm{spec}},$$

where we know, by Lemmas B.2 and B.3, that $\|U_{j,\mathcal{R};T}\|_2^2 = O(2^j T^{-1})$ and $\|\Sigma_T\|_{\mathrm{spec}} \leq \|c_X\|_{1,\infty}$. Moreover, we can write

$$\|\bar{\Sigma}_T - \Sigma_T\|_{\mathrm{spec}}$$
$$\leq \sum_{u=-\infty}^{\infty} \sup_s (\sigma_{s,s+u} - \bar{\sigma}_{s,s+u})$$
(B.14)
$$= \sum_{u=-\infty}^{\infty} \sup_s \sum_{\ell=-\infty}^{-1} \sum_{n=-\infty}^{\infty} (w_{\ell n;T}^2 - Q_{\ell,\mathcal{R}_T(s)})$$
$$\qquad \times \psi_{\ell n}(s)\psi_{\ell n}(s+u) + \mathrm{R}_1 + \mathrm{R}_2,$$

where

$$\mathrm{R}_1 = \sum_{u=-\infty}^{\infty} \sup_s \sum_{\ell=-\infty}^{-1} Q_{\ell,\mathcal{R}_T(s)} \Psi_\ell(u) \mathbb{I}_{|u|>M_T},$$
$$\mathrm{R}_2 = \sum_{u=-\infty}^{\infty} \sup_s \sum_{\ell=-\infty}^{-\log_2(T)-1} Q_{\ell,\mathcal{R}_T(s)} \Psi_\ell(u) \mathbb{I}_{|u|<M_T}.$$

As

$$\sum_{u=-\infty}^{\infty} \sup_s \sum_{\ell=-\infty}^{-1} Q_{\ell,\mathcal{R}_T(s)} \Psi_\ell(u) = \sum_{u=-\infty}^{\infty} \sup_s |\mathcal{R}_T|^{-1} \int_{\mathcal{R}_T(s)} dz\, c_X(z,u),$$

the rate of $R_1$ is $o_T(2^{-J_T})$, by Assumption 6. Next, using $|\Psi_\ell(u)| \leq 1$ uniformly in $\ell < 0$, we get

$$|R_2| \leq \sum_{u=-\infty}^{\infty} \sup_s |\mathcal{R}_T|^{-1} \int_{\mathcal{R}_T(s)} dz \sum_{\ell=-\infty}^{-\log_2(T)-1} S_\ell(z) \mathbb{I}_{|u|<M_T}$$

$$\leq 2M_T \sum_{\ell=-\infty}^{-\log_2(T)-1} \sup_z S_\ell(z) = O(M_T/T),$$

using Assumption 4. Assumption 5 on the rate of the truncating sequence $M_T$ implies $|R_2| = o_T(2^{-J_T})$. The main term of (B.14) is bounded by

(B.15)
$$\sum_{u=-\infty}^{\infty} \sup_s \sum_{\ell=-\infty}^{-1} \sum_{n=-\infty}^{\infty} |\mathcal{R}_T|^{-1} \int_{\mathcal{R}_T(s)} dz |w_{\ell n;T}^2 - S_\ell(z)|$$
$$\times |\psi_{\ell n}(s)\psi_{\ell n}(s+u)|.$$

By Definition 1, we can write

$$|w_{\ell n;T}^2 - S_\ell(z)| \leq \frac{\overline{C}C_\ell}{T} + \left|S_\ell\left(\frac{n}{T}\right) - S_\ell\left(\frac{n-s}{T} + z\right)\right|$$
$$+ \left|S_\ell(z) - S_\ell\left(\frac{n-s}{T} + z\right)\right|,$$

which, when substituted into (B.15), leads to three terms. By (B.6) and (2.4), the first term is $O(T^{-1})$. For the second term, with a change of variable $z$ to $z + s/T$, we get

$$\sum_{u=-\infty}^{\infty} \sup_s \sum_{\ell=-\infty}^{-1} \sum_{n=-\infty}^{\infty} |\mathcal{R}_T|^{-1} \int_{\mathcal{R}_T(0)} dz \left|S_\ell\left(\frac{n}{T}\right) - S_\ell\left(\frac{n}{T} + z\right)\right|$$
$$\times |\psi_{\ell n}(s)\psi_{\ell n}(s+u)|,$$

where $\mathcal{R}_T(0)$ denotes the interval $\mathcal{R}_T(s)$ shifted by $-s$. If we use the fact that $|\psi_{\ell n}(s)|$ is uniformly bounded and $\sum_{u=-\infty}^{\infty} |\psi_{\ell n}(s+u)| = O(\mathcal{L}_\ell)$, the second term is then bounded (up to a multiplicative constant) by

$$|\mathcal{R}_T|^{-1} \sum_{\ell=-\infty}^{-1} \mathcal{L}_\ell \int_{\mathcal{R}_T(0)} dz \sum_{n=-\infty}^{\infty} \left|S_\ell\left(\frac{n}{T}\right) - S_\ell\left(\frac{n}{T} + z\right)\right|$$
$$\leq |\mathcal{R}_T|^{-1} \sum_{\ell=-\infty}^{-1} \mathcal{L}_\ell \int_{\mathcal{R}_T(0)} dz |z| \operatorname{TV}(S_\ell)$$
$$= O(|\mathcal{R}_T|) \sum_{\ell=-\infty}^{-1} \mathcal{L}_\ell L_\ell = O(|\mathcal{R}_T|),$$

by assumptions (2.3) and (2.4). The third term is

$$\sum_{u=-\infty}^{\infty} \sup_s \sum_{\ell=-\infty}^{-1} \sum_{n=-\infty}^{\infty} |\mathcal{R}_T|^{-1} \int_{\mathcal{R}_T(s)} dz \left| S_\ell(z) - S_\ell\left(\frac{n-s}{T} + z\right) \right|$$
$$\times |\psi_{\ell n}(s)\psi_{\ell n}(s+u)|.$$

If $s_0$ denotes the infimum of $\mathcal{R}_T(s)$, we decompose the integral as follows:

$$\sum_{u=-\infty}^{\infty} \sup_s \sum_{\ell=-\infty}^{-1} \sum_{n=-\infty}^{\infty} |\mathcal{R}_T|^{-1} \sum_{k=0}^{|\mathcal{R}_T T|-1} \int_{s_0+k/T}^{s_0+(k+1)/T} dz \left| S_\ell(z) - S_\ell\left(\frac{n-s}{T} + z\right) \right|$$
$$\times |\psi_{\ell n}(s)\psi_{\ell n}(s+u)|,$$

which can be rewritten, with the change of variable $y := z - s_0 - k/T$, as

$$\sum_{u=-\infty}^{\infty} \sup_s \sum_{\ell=-\infty}^{-1} \sum_{n=-\infty}^{\infty} |\mathcal{R}_T|^{-1}$$
$$\times \sum_{k=0}^{|\mathcal{R}_T T|-1} \int_0^{1/T} dy \left| S_\ell\left(y + s_0 + \frac{k}{T}\right) - S_\ell\left(y + s_0 + \frac{n-s+k}{T}\right) \right|$$
$$\times |\psi_{\ell n}(s)\psi_{\ell n}(s+u)|.$$

Assumption (2.3) for the sum over $k$ with (B.5) leads to the bound

$$\sum_{u=-\infty}^{\infty} \sup_s \sum_{\ell=-\infty}^{-1} L_\ell \sum_{n=-\infty}^{\infty} |\mathcal{R}_T T|^{-1} |n-s| |\psi_{\ell n}(s)\psi_{\ell n}(s+u)|.$$

The compact support of $\psi_{\ell n}(s)$ implies $|n - s| < \mathcal{L}_\ell$. Therefore, (B.6), (2.3) and (2.4) imply that this last term is $O(|\mathcal{R}_T T|^{-1})$. Finally, we summarize all the rates of convergence for the deterministic term as follows:

$$2^{-j}T \cdot (\bar{\sigma}^2_{j,\mathcal{R},T} - \sigma^2_{j,\mathcal{R},T})$$
$$= O(T^{-1} + |\mathcal{R}_T| + |\mathcal{R}_T T|^{-1}) + |\mathrm{R}_1| + |\mathrm{R}_2|$$
$$= O(T^{-1} + |\mathcal{R}_T| + |\mathcal{R}_T T|^{-1}) + o_T(2^{-J_T}) + o_T(2^{-J_T})$$
$$= o_T(2^{-J_T}),$$

by Assumption 5.

Let us now turn to the stochastic term $|\tilde{\sigma}^2_{j,\mathcal{R},T} - \bar{\sigma}^2_{j,\mathcal{R},T}|$. Lemma B.5 implies the existence of a random set $\mathcal{A}$ which does not depend on $j$ and such that $\Pr(\mathcal{A}) \geq 1 - o_T(1)$ and $|Q_{j,\mathcal{R}_T(s);T} - Q_{j,\mathcal{R}_T(s)}| \leq (k_T/|\mathcal{R}_T|)\sqrt{(C^2+c^2)/T}$ almost surely on $\mathcal{A}$, for all $T > T_0$ and $j = -1, \ldots, -J_T$. We can write

$$|\tilde{\sigma}^2_{j,\mathcal{R},T} - \bar{\sigma}^2_{j,\mathcal{R},T}|$$

$$\leq 2 \sum_{h,t=0}^{T-1} \left| \sum_{s,u=0}^{T-1} U_{ts}^{(j)} U_{tu}^{(j)} \right.$$

(B.16)
$$\times \sum_{\ell,m=-\log_2 T}^{-1} (Q_{\ell,\mathcal{R}_T(s);T} Q_{m,\mathcal{R}_T(u);T}$$

$$\left. - Q_{\ell,\mathcal{R}_T(s)} Q_{m,\mathcal{R}_T(u)}) \Psi_\ell(s-h)\Psi_m(u-h) \right|$$

$$\times \mathbb{I}_{|s-h|\leq M_T}\mathbb{I}_{|u-h|\leq M_T}$$

almost surely on $\mathcal{A}$. Using the decomposition

$$Q_{\ell,\mathcal{R}_T(s);T} Q_{m,\mathcal{R}_T(u);T} - Q_{\ell,\mathcal{R}_T(s)} Q_{m,\mathcal{R}_T(u)}$$
$$= (Q_{m,\mathcal{R}_T(u);T} - Q_{m,\mathcal{R}_T(u)}) Q_{\ell,\mathcal{R}_T(s)}$$
$$+ (Q_{\ell,\mathcal{R}_T(s);T} - Q_{\ell,\mathcal{R}_T(s)}) Q_{m,\mathcal{R}_T(u)}$$
$$+ (Q_{\ell,\mathcal{R}_T(s);T} - Q_{\ell,\mathcal{R}_T(s)})(Q_{m,\mathcal{R}_T(u);T} - Q_{m,\mathcal{R}_T(u)}),$$

we get three terms in the right-hand side of (B.16). On $\mathcal{A}$, the first of these terms is bounded as follows (the other terms are bounded similarly):

$$2 \sum_{h,t,s,u} \left| U_{ts}^{(j)} U_{tu}^{(j)} \sum_m (Q_{m,\mathcal{R}_T(u);T} - Q_{m,\mathcal{R}_T(u)}) \right.$$

$$\left. \times \Psi_m(u-h) \sum_\ell Q_{\ell,\mathcal{R}_T(s)} \Psi_\ell(s-h) \right| \mathbb{I}_{|s-u|\leq 2M_T}$$

$$\leq 2\sqrt{1+c^2} \frac{k_T \log_2 T}{|\mathcal{R}_T|\sqrt{T}} \sum_{h,t,s,u} |U_{ts}^{(j)} U_{tu}^{(j)}|$$

$$\times \left| \sum_\ell Q_{\ell,\mathcal{R}_T(s)} \Psi_\ell(s-h) \right| \mathbb{I}_{|s-u|\leq 2M_T}$$

$$\leq 2\sqrt{1+c^2} \frac{k_T \log_2 T}{|\mathcal{R}_T|\sqrt{T}} \sum_{t,s,u} |U_{ts}^{(j)} U_{tu}^{(j)}| \mathbb{I}_{|s-u|\leq 2M_T}$$

$$\times \sum_h \sup_z \left| \sum_\ell S_\ell(z) \Psi_\ell(h) \right|$$

$$= O(2^j M_T k_T |\mathcal{R}_T T|^{-1} T^{-1/2} \log_2^3 T) \quad \text{a.s. on } \mathcal{A},$$

using Assumption 1 and Lemma B.4. The result then follows from Assumption 5. $\square$

**B.6. Proof of Theorem 1.** By Proposition 4 and for $T$ large enough, there exists a random set $\mathcal{A}$ such that $1 - \Pr(\mathcal{A}) = o_T(1)$ and (3.9) holds on $\mathcal{A}$. Then, if $\mathcal{A}^c$ denotes the complementary random set of $\mathcal{A}$, we can write

$$\Pr(|Q_{j,\mathcal{R};T} - Q_{j,\mathcal{R}}| > 2\tilde{\sigma}_{j,\mathcal{R},T}\eta)$$
$$= \Pr(|Q_{j,\mathcal{R};T} - Q_{j,\mathcal{R}}| > 2\tilde{\sigma}_{j,\mathcal{R},T}\eta|\mathcal{A})\Pr(\mathcal{A})$$
$$+ \Pr(|Q_{j,\mathcal{R};T} - Q_{j,\mathcal{R}}| > 2\tilde{\sigma}_{j,\mathcal{R},T}\eta|\mathcal{A}^c)(1 - \Pr(\mathcal{A})).$$

The second term of this sum is $o_T(1)$, by Proposition 4. To bound the first term, we observe that Proposition 4 implies $\tilde{\sigma}^2_{j,\mathcal{R},T} \geq \sigma^2_{j,\mathcal{R},T} - \varphi_T$ on $\mathcal{A}$ with $\varphi_T = o_T(2^{j-J_T}T^{-1})$. Together with Proposition 2, this implies

$$(B.17) \qquad \frac{\tilde{\sigma}^2_{j,\mathcal{R},T}}{\sigma^2_{j,\mathcal{R},T}} \geq 1 - \frac{\varphi_T}{\sigma^2_{j,\mathcal{R},T}} = 1 - o_T(1) \to 1$$

for all $j = -1, \ldots, -J_T$, as $T$ tends to infinity. We can then write

$$\Pr(|Q_{j,\mathcal{R};T} - Q_{j,\mathcal{R}}| > 2\tilde{\sigma}_{j,\mathcal{R},T}\eta)$$
$$\leq \Pr\left(|Q_{j,\mathcal{R};T} - Q_{j,\mathcal{R}}| > 2\sigma_{j,\mathcal{R},T}\eta\sqrt{1 - \frac{\varphi_T}{\sigma^2_{j,\mathcal{R},T}}}\Big|\mathcal{A}\right) + o_T(1)$$

and Proposition 3 leads to the result with $\gamma_T = \varphi_T/\sigma^2_{j,\mathcal{R};T}$.

**B.7. Proof of Proposition 5.** Let $\mathcal{U}$ be a segment of $\wp(\mathcal{R})$. Consider the a.s. inequality

$$|Q_{j,\mathcal{R};T} - Q_{j,\mathcal{U};T}| \leq |Q_{j,\mathcal{R};T} - Q_{j,\mathcal{R}}| + |Q_{j,\mathcal{U};T} - Q_{j,\mathcal{U}}| + \Delta_j(\mathcal{R},\mathcal{U}),$$

where $\Delta_j(\mathcal{R},\mathcal{U})$ is defined in (4.1). In the regular case, $\Delta_j(\mathcal{R},\mathcal{U}) \leq b(\mathcal{U}) + b(\mathcal{R}) \leq C_j(\sigma_{j,\mathcal{U},T} + \sigma_{j,\mathcal{R},T})k_T$. Consequently, in the regular case,

$\Pr(\mathcal{R} \text{ is rejected})$

$$\leq \sum_{\mathcal{U} \in \wp(\mathcal{R})} \Pr\{|Q_{j,\mathcal{U};T} - Q_{j,\mathcal{R};T}| > 2(\eta\sigma_{j,\mathcal{U},T} + \eta\sigma_{j,\mathcal{R},T})k_T\}$$

$$\leq \sum_{\mathcal{U} \in \wp(\mathcal{R})} \Pr(|Q_{j,\mathcal{R};T} - Q_{j,\mathcal{R}}| > -C_j\sigma_{j,\mathcal{R},T}k_T + 2\eta\sigma_{j,\mathcal{R},T}k_T)$$

$$+ \sum_{\mathcal{U} \in \wp(\mathcal{R})} \Pr(|Q_{j,\mathcal{U};T} - Q_{j,\mathcal{U}}| > -C_j\sigma_{j,\mathcal{U},T}k_T + 2\eta\sigma_{j,\mathcal{U},T}k_T).$$

Proposition 3 implies

$\Pr(\mathcal{R} \text{ is rejected})$

$$\leq (\sharp\wp(\mathcal{R}))c_0 \exp\left\{-\frac{1}{16} \cdot \eta_T^2 \Big/ \left[1 + \frac{2\eta_T L_j}{|\mathcal{R}T|\sigma_{j,\mathcal{R},T}}\right.\right.$$

$$+ \frac{2^{j/2}\eta_T(K_2\|c_X\|_{1,\infty}+K_3)}{\sigma_{j,\mathcal{R},T}|\mathcal{R}|\sqrt{T}}\Bigg]\Bigg\}$$

$$+ c_0 \sum_{\mathcal{U}\in\wp(\mathcal{R})} \exp\Bigg\{-\frac{1}{16}\cdot\eta_T^2\Big/\Bigg[1+\frac{2\eta_T L_j}{|\mathcal{U}T|\sigma_{j,\mathcal{U},T}}$$

$$+ \frac{2^{j/2}\eta_T(K_2\|c_X\|_{1,\infty}+K_3)}{\sigma_{j,\mathcal{U},T}|\mathcal{U}|\sqrt{T}}\Bigg]\Bigg\},$$

with $\eta_T := 2\eta k_T - C_j k_T = k_T 2^{-j/2}(5(2\alpha+p) - \sqrt{\alpha+p})$.

Proposition 2 leads to $\sigma_{j,\mathcal{R};T}^{-1} \leq C^{-1}2^{-j/2}\sqrt{|\mathcal{R}T|}$ and similarly for $\sigma_{j,\mathcal{U},T}^{-1}$. As $\delta \leq |\mathcal{U}| \leq |\mathcal{R}| \leq 1$, we consider the dominant terms in the sum and can write, for $T$ large enough and with $2^{-j/2}L_j \leq \rho\mathcal{L}_{-1}$,

$\Pr(\mathcal{R}$ is rejected)

$$\leq 2c_0(\sharp\wp(\mathcal{R}))\exp\Bigg\{-\frac{1}{16}\cdot\eta_T^2\Big/\Bigg[1+\frac{2\eta_T\rho\mathcal{L}_{-1}}{\sqrt{K_1\delta T}}$$

$$+ \frac{\eta_T(K_2\|c_X\|_{1,\infty}+K_3)}{\sqrt{2^j K_1\delta}}\Bigg]\Bigg\}.$$

Replacing $\eta_T$, using $2\alpha + p \geq \sqrt{\alpha+p}$ and $k_T \sim \log_2 T$, the asymptotic order of this bound is

$$(\sharp\wp(\mathcal{R}))O(T^{-(\sqrt{\delta K_1}/(K_2\|c_X\|_{1,\infty}+K_3))(\alpha+p/2)})$$

and the result follows for $T$ large enough by Assumption 7(2).

**B.8. Proof of Theorem 2.** For reader's convenience, we first state two technical lemmas. The first lemma is a consequence of Rosenthal's inequality [see, e.g., Härdle, Kerkyacharian, Picard and Tsybakov (1998)].

LEMMA B.6. *Let $Y \sim \mathcal{N}(0,\sigma^2)$ with $\sigma^2 > 0$. Then, $\mathrm{E}|Y|^p \leq C(p)\sigma^p$, where $C(p)$ is a function of $p$ only.*

LEMMA B.7. *Let $\underline{Z}_T = (Z_1,\ldots,Z_T)'$ be a vector of i.i.d. Gaussian random variables with zero mean and $\mathrm{Var}\, Z_1 = 1$. If $M_{j,\mathcal{R};T}$ is the matrix (B.9), $v$ is a positive constant and $p \geq 2$, then there exists $T_0$ such that*

$$\mathrm{E}(\underline{Z}_T' M_{j,\mathcal{R};T}\underline{Z}_T - \mathrm{tr}\, M_{j,\mathcal{R};T} + vk_T T^{-1/2})^p$$

$$\leq C(\kappa,\|c_X\|_{1,\infty},p)T^{-p/2}(2^{1+j/2}|\mathcal{R}|^{-1}+vk_T)^p$$

*for all $T \geq T_0$.*

PROOF. First, we write

(B.18)
$$E(\underline{Z}'_T M_{j,\mathcal{R};T}\underline{Z}_T - \text{tr}\, M_{j,\mathcal{R};T} + vk_T T^{-1/2})^p$$
$$= \sum_{r=0}^{p} \binom{p}{r} E(\underline{Z}'_T M_{j,\mathcal{R};T}\underline{Z}_T - \text{tr}\, M_{j,\mathcal{R};T})^r v^{p-r} k_T^{p-r} T^{-(p-r)/2}.$$

Due to the relationship between the centered moments of a random variable and its cumulants, we can write

$$E(\underline{Z}'_T M_{j,\mathcal{R};T}\underline{Z}_T - \text{tr}\, M_{j,\mathcal{R};T})^r$$
$$= \sum_{m=0}^{r} \sum C(p_1,\ldots,p_m,m,\pi_1,\ldots,\pi_m,r)\kappa_{p_1}^{\pi_1}\ldots\kappa_{p_m}^{\pi_m},$$

where the second sum is over $p_1,\ldots,p_m,\pi_1,\ldots,\pi_m$ in $\{1,\ldots,r\}$ such that $\sum_{i=1}^{m} p_i \pi_i = r$, $\kappa_{p_i}$ is the $p_i$th cumulant of $\underline{Z}'_T M_{j,\mathcal{R};T}\underline{Z}_T$ and $C$ denotes a generic constant in this proof. From Lemma B.1, (B.13) and Proposition 2, $\kappa_{p_i} \leq 2^{p_i} \times (p_i-1)! K_2^{p_i} \|c_X\|_{1,\infty}^{p_i} 2^{jp_i/2} |\mathcal{R}|^{-p_i} T^{-p_i/2}$ and, consequently, $E(\underline{Z}'_T M_{j,\mathcal{R};T}\underline{Z}_T - \text{tr}\, M_{j,\mathcal{R};T})^r \leq C(\kappa, \|c_X\|_{1,\infty}, r)\, 2^{r(1+j/2)}|\mathcal{R}|^{-r}T^{-r/2}$. Using this inequality in (B.18) leads to the result. $\square$

PROOF OF THEOREM 2. In this proof, $C$ denotes a generic constant. Let $\tilde{\mathcal{R}}$ be the interval selected by the estimation procedure. We consider two cases, $|\tilde{\mathcal{R}}| < |\mathcal{R}|$ and $|\tilde{\mathcal{R}}| \geq |\mathcal{R}|$, and split the expectation into two parts as follows:

$$E|\tilde{S}_j(z_0) - S_j(z_0)|^p$$
$$= E|\tilde{S}_j(z_0) - S_j(z_0)|^p 1_{|\tilde{\mathcal{R}}|<|\mathcal{R}|} + E|\tilde{S}_j(z_0) - S_j(z_0)|^p 1_{|\tilde{\mathcal{R}}|\geq|\mathcal{R}|}.$$

*First term* ($|\tilde{\mathcal{R}}| < |\mathcal{R}|$). In the first case, we make use of the inequality $|a-b|^p \leq 2^{p-1}|a|^p + 2^{p-1}|b|^p$ and write

$$E|\tilde{S}_j(z_0) - S_j(z_0)|^p 1_{|\tilde{\mathcal{R}}|<|\mathcal{R}|}$$
$$\leq 2^{p-1} E|S_j(z_0) - Q_{j,\tilde{\mathcal{R}}}|^p 1_{|\tilde{\mathcal{R}}|<|\mathcal{R}|} + 2^{p-1} E|Q_{j,\tilde{\mathcal{R}};T} - Q_{j,\tilde{\mathcal{R}}}|^p 1_{|\tilde{\mathcal{R}}|<|\mathcal{R}|}.$$

As $|\tilde{\mathcal{R}}| < |\mathcal{R}|$, the evolutionary wavelet spectrum is homogeneous over $\mathcal{R}$ and $\tilde{\mathcal{R}}$, and property (4.4) holds for $\tilde{\mathcal{R}}$. Then, using Proposition 2 on the variance and the first point of Assumption 7, the first term of the right-hand side is bounded as follows:

(B.19)
$$2^{p-1} E|S_j(z_0) - Q_{j,\tilde{\mathcal{R}}}|^p 1_{|\tilde{\mathcal{R}}|<|\mathcal{R}|}$$
$$\leq 2^{p-1} E(C_j \sigma_{j,\tilde{\mathcal{R}},T} k_T)^p$$
$$\leq 2^{p-1} C_j^p k_T^p 2^{jp/2}(T\delta^2)^{-p/2}(1+c^2)^{p/2}$$
$$= 2^{p-1}(\alpha+p)^{p/2} k_T^p T^{-p/2}\delta^{-p}(1+c^2)^{p/2},$$

by the definition of $C_j$ [see equation (4.5)]. Now, if we let $G_T = \underline{Z}'_T M_{j,\tilde{\mathcal{R}};T} \underline{Z}_T + |\tilde{\mathcal{R}}T|^{-1} \sum_{k \in \tilde{\mathcal{R}}T} z_{j,k;T} - \operatorname{tr} M_{j,\tilde{\mathcal{R}};T}$, then the second term may be written

$$2^{p-1} \operatorname{E} |G_T + \operatorname{bias}_T|^p 1_{|\tilde{\mathcal{R}}|<|\mathcal{R}|} \leq 2^{2p-2} \{\operatorname{E}(|G_T|^p 1_{|\tilde{\mathcal{R}}|<|\mathcal{R}|}) + |\operatorname{bias}_T|^p\},$$

where, using Proposition 2 for $T$ large enough,

(B.20) $$|\operatorname{bias}_T|^p \leq C_p 2^{jp/2} (\delta T)^{-p/2},$$

with a constant $C_p$ depending only on $p$. Finally, we now show that $\operatorname{E}|G_T|^p$ is uniformly bounded in $T$. Using $\delta < |\tilde{\mathcal{R}}| < |\mathcal{R}|$, we first note that Propositions 2 and B.1 imply

(B.21) $$\Pr\left(|G_T| > \frac{\lambda}{\delta}\sqrt{(C^2+c^2)\frac{2^j}{T}}\right) \leq 3\exp\left(-\frac{1}{4} \cdot \frac{\lambda^2}{1 + 2\lambda\tau_\infty\sqrt{|\mathcal{R}T|/2^j}}\right),$$

where $\tau_\infty \leq 2^{(j-1)/2} c/(\delta\sqrt{T})$ by (B.13). We now truncate the integral $\operatorname{E}|G_T|^p = \int_0^\infty dx \Pr(|G_T|^p \geq x)$ at the point $\mu_T^{p/2}$, which is such that $\mu_T = 2^j(C^2 + c^2)/(\delta^2 T)$. With the change of variable $x = y^p \mu_T^{p/2}$, this leads to

$$\operatorname{E}|G_T|^p \leq \mu_T^{p/2} + p\mu_T^{p/2} \int_1^\infty dy \, y^{p-1} \Pr(|G_T| > y\mu_T^{1/2})$$

$$\leq \mu_T^{p/2} + p\mu_T^{p/2} \int_1^\infty dy \, y^{p-1} \exp\left(-\frac{1}{2} \cdot \frac{y^2}{1+2y\tau_\infty\sqrt{|\mathcal{R}T|/2^j}}\right).$$

To compute the integral, we note that $1 \leq y$ and evaluate $\int_1^\infty dy\, y^{p-1} \exp(-\alpha_T y)$. This leads to the bound

$$\operatorname{E}|G_T|^p \leq \mu_T^{p/2} + ep\mu_T^{p/2}\left(2 + 4\tau_\infty\sqrt{\frac{|\mathcal{R}T|}{2^j}}\right)^p \leq C_p \delta^{-p} T^{-p/2}.$$

In conclusion, in the first case, we get the bound $\operatorname{E}|\tilde{S}_j(z_0) - S_j(z_0)|^p 1_{|\tilde{\mathcal{R}}|<|\mathcal{R}|} \leq C_p \delta^{-p} T^{-p/2} k_T^p$ from (B.19) and (B.20).

*Second term* ($|\tilde{\mathcal{R}}| \geq |\mathcal{R}|$). We now consider the second case. Select a subinterval $\mathcal{U}$ in $\wp(\tilde{\mathcal{R}})$ included in $\mathcal{R}$ and containing $z_0$. Then, consider the decomposition

$$\operatorname{E}|\tilde{S}_j(z_0) - S_j(z_0)|^p 1_{|\tilde{\mathcal{R}}|\geq|\mathcal{R}|}$$
$$\leq \operatorname{E}\{|Q_{j,\mathcal{U}} - S_j(z_0)| + |Q_{j,\mathcal{U};T} - Q_{j,\mathcal{U}}| + |Q_{j,\tilde{\mathcal{R}};T} - Q_{j,\mathcal{U};T}|\}^p.$$

As the wavelet spectrum is regular on $\mathcal{U} \subset \mathcal{R}$, the term $|Q_{j,\mathcal{U}} - S_j(z_0)|$ is bounded by $C_j \sigma_{j,\mathcal{U},T} k_T$. On the other hand, using Proposition 2, $|Q_{j,\mathcal{U};T} - Q_{j,\mathcal{U}}| = |Q_{j,\mathcal{U};T} - \operatorname{tr} M_{j,\mathcal{U};T}| + R_T$ with $R_T = O(2^{j/2} T^{-1/2})$. Moreover, as $\tilde{\mathcal{R}}$

is selected by the estimation procedure, it holds that $|Q_{j,\tilde{\mathcal{R}};T} - Q_{j,\mathcal{U};T}| \leq 2(\eta\sigma_{j,\tilde{\mathcal{R}},T} + \lambda\sigma_{j,\mathcal{U},T})k_T$ almost surely. With $2\alpha + p \geq \sqrt{\alpha + p}$, we can write

$$C_j \sigma_{j,\mathcal{U},T} k_T + 2(\eta \sigma_{j,\tilde{\mathcal{R}},T} + \lambda \sigma_{j,\mathcal{U},T}) k_T$$
$$\leq 11\sqrt{2}(2\alpha + p) k_T (1 + c^2) T^{-1/2} \delta^{-1},$$

using $|\tilde{\mathcal{R}}| \geq |\mathcal{U}| \geq \delta$. Then, Lemmas B.6 and B.7 prove the existence of a constant $c_5$ depending on $\kappa, \nu, p, K_2$ and on $\|c_X\|_{1,\infty}$, such that, for $T \geq T_0$,

$$\mathrm{E}\{|Q_{j,\mathcal{U};T} - \mathrm{tr}\, M_{j,\mathcal{U};T}| + R_T + C_j \sigma_{j,\mathcal{U},T} k_T + 2(\eta \sigma_{j,\tilde{\mathcal{R}},T} + \lambda \sigma_{j,\mathcal{U},T}) k_T\}^p$$
$$\leq C_p \delta^{-p} T^{-p/2} (2^{j/2} |\mathcal{U}|^{-1} + k_T)^p + C_p 2^{jp/2} |\mathcal{U}T|^{-p/2}$$

and the result follows using $|\mathcal{U}| \geq \delta$. □

**B.9. Proof of Proposition 6.** We first prove the following lemma, stating an exponential inequality for quadratic forms of Gaussian random variables.

LEMMA B.8. *Let $\underline{Z}_T = (Z_1, \ldots, Z_T)'$ be a vector of i.i.d. Gaussian random variables with zero mean and $\mathrm{Var}\, Z_1 = 1$. If $M_T$ is a $T \times T$ symmetric and positive definite matrix, then*

$$\Pr(\underline{Z}_T' M_T \underline{Z}_T \leq \eta) \leq \exp\left(-\frac{(\eta - \mathrm{tr}\, M_T)^2}{4\|M_T\|_2^2}\right),$$

*provided that $\eta \leq \mathrm{tr}\, M_T$.*

PROOF. By assumption on the matrix $M_T$, the decomposition $M_T = O_T' \Lambda_T \times O_T$ holds with a diagonal $T \times T$ matrix $\Lambda_T$ and an orthonormal matrix $O_T$. If we let $\underline{Y}_T = O_T' \underline{Z}_T$, then $\underline{Y}_T$ is a vector of i.i.d. Gaussian random variables with zero mean and $\mathrm{Var}\, Y_1 = 1$. We can write $\underline{Z}_T' M_T \underline{Z}_T = \underline{Y}_T' \Lambda_T \underline{Y}_T = \sum_{i=1}^T \lambda_i Y_i^2$, with $\lambda_i > 0$. Moreover, $\mathrm{tr}\, M_T = \mathrm{tr}\, \Lambda_T$, $\mathrm{tr}\, \Lambda_T^2 = \mathrm{tr}\, M_T^2 = \|M_T\|_2^2$ and $\|M_T\|_{\mathrm{spec}} = \max\{\lambda_1, \ldots, \lambda_T\}$. The Chernoff inequality [Ross (1998)] on $\underline{Y}_T$ leads to

$$\Pr(\underline{Z}_T' M_T \underline{Z}_T \leq \eta) = \Pr(\underline{Y}_T' \Lambda_T \underline{Y}_T \leq \eta)$$
$$\leq \exp\left\{\inf_{t<0}(-t\eta + \log \mathrm{E}\exp(t\underline{Y}_T' \Lambda_T \underline{Y}_T))\right\}$$
$$= \exp\left\{\inf_{t<0}\left(-t\eta + \sum_{i=1}^T \log \mathrm{E}\exp(\lambda_i t Y_i^2)\right)\right\}$$

and, using the fact that $\log \mathrm{E}\exp(\alpha_i Y_i^2) = -\frac{1}{2}\log(1 - 2\alpha_i) \leq \alpha_i + \alpha_i^2$ holds for $\alpha_i \leq 0$, we get

$$\Pr(\underline{Z}_T' M_T \underline{Z}_T \leq \eta) \leq \exp\left\{\inf_{t<0}(-t\eta + t\,\mathrm{tr}\, \Lambda_T + t^2 \mathrm{tr}\, \Lambda_T^2)\right\}.$$

The result follows by taking $t = (\eta - \operatorname{tr} \Lambda_T)/(2 \operatorname{tr} \Lambda_T^2)$. □

Lemma B.8 is not directly applicable to the quadratic form $Q_{j,\mathcal{R};T} = \underline{Z}_T' M_{j,\mathcal{R};T} \underline{Z}_T$ because the matrix $M_{j,\mathcal{R};T}$ is not positive definite in general. In the next lemma, we show how this matrix can be approximated by the matrix $M^\star_{j,\mathcal{R};T}$, defined as

$$M^\star_{j,\mathcal{R};T} = \Sigma_T^{1/2\prime} U^\star_{j,\mathcal{R};T} \Sigma_T^{1/2},$$

where the entry $(s,t)$ of the matrix $U^\star_{j,\mathcal{R};T}$ is given by

$$u^\star_{st} = 2\gamma_0 |\mathcal{R}T|^{-1} \sum_{\ell=-\log_2 T}^{-1} 2^{\ell/2} \Psi_\ell(s-t),$$

with $\gamma_0 \geq \sup_{j<0} \sup_{\ell<0} 2^{-\ell/2} |A_{j\ell}^{-1}| > 0$. The matrix $M^\star_{j,\mathcal{R};T}$ is clearly symmetric. It is also positive definite because $U^\star_{j,\mathcal{R};T}$ is positive definite: for all sequences $\underline{x} = (x_0, \ldots, x_{T-1})'$ of $\ell^2$, the quadratic form

$$x' U^\star_{j,\mathcal{R};T} x = \gamma_0 |\mathcal{R}T|^{-1} \sum_{\ell=-\log_2 T}^{-1} 2^{\ell/2} \sum_s \left( \sum_k x_s \psi_{\ell k}(s) \right)^2$$

is strictly positive.

LEMMA B.9. *Suppose that Assumptions 1–4 hold true. Define $\gamma_1$ such that*

$$0 < \gamma_1 \leq \gamma_0 \inf_{m<0} \sum_{\ell=-\log_2 T}^{-1} 2^{\ell/2} A_{m\ell}.$$

*The following properties hold true for $T$ sufficiently large:*

(B.22) $\quad \gamma_1 |\mathcal{R}|^{-1} \varepsilon \leq \operatorname{tr}(M^\star_{j,\mathcal{R};T} - M_{j,\mathcal{R};T}) \leq 6 \|c_{X,T}\|_{1,\infty} \gamma_0 |\mathcal{R}|^{-1},$

*where $\varepsilon$ is defined in Assumption 2,*

(B.23)
$$\begin{aligned}
&\|M^\star_{j,\mathcal{R};T} - M_{j,\mathcal{R};T}\|^2_{\mathrm{spec}} \\
&\leq \|M^\star_{j,\mathcal{R};T} - M_{j,\mathcal{R};T}\|^2_2 \\
&\leq 8\mathcal{L}_{-1} \gamma_0^2 |\mathcal{R}|^{-2} \|c_X\|^2_{1,\infty} T^{-1} \log_2^2(T) + O(T^{-1})
\end{aligned}$$

*and, if $\underline{Z}_T = (Z_1, \ldots, Z_T)'$ is a vector of i.i.d. Gaussian random variables with zero mean and $\operatorname{Var} Z_1 = 1$, then*

(B.24) $\quad \Pr(\underline{Z}_T'(M^\star_{j,\mathcal{R};T} - M_{j,\mathcal{R};T})\underline{Z}_T > \lambda_T) = O\left(\exp\left\{-\frac{\sqrt{T} \operatorname{tr} M_{j,\mathcal{R};T}}{\log_2^2 T}\right\}\right),$

*where $\lambda_T = \operatorname{tr} M^\star_{j,\mathcal{R};T} - \operatorname{tr} M_{j,\mathcal{R};T} + \operatorname{tr} M_{j,\mathcal{R};T} \log_2^{-1} T$.*

PROOF. 1. We prove (B.22). Write $\operatorname{tr}(M^\star_{j,\mathcal{R};T} - M_{j,\mathcal{R};T}) = \operatorname{tr}(M^\star_{j,\mathcal{R};T}) - \operatorname{tr}(M_{j,\mathcal{R};T})$, where the second term is $\mathrm{E}(\underline{Z}'_T M_{j,\mathcal{R};T} \underline{Z}_T) = Q_{j,\mathcal{R}} + O(2^{j/2} T^{-1/2})$, from Lemma B.1 and Proposition 2. Moreover,

$$\operatorname{tr}(M^\star_{j,\mathcal{R};T})$$
$$= \operatorname{tr}(\Sigma'_T U^\star_{j,\mathcal{R};T})$$

$$(B.25) \qquad = 2\gamma_0 |\mathcal{R}T|^{-1} \sum_{s,u=-\infty}^{\infty} c_{X,T}\left(\frac{s}{T}, u\right) \sum_{\ell=-\log_2 T}^{-1} 2^{\ell/2} \Psi_\ell(u)$$

$$(B.26) \qquad = 2\gamma_0 |\mathcal{R}T|^{-1} \sum_{s,u=-\infty}^{\infty} c_X\left(\frac{s}{T}, u\right) \sum_{\ell=-\log_2 T}^{-1} 2^{\ell/2} \Psi_\ell(u) + \operatorname{Rest}_T.$$

We now derive a bound for $\operatorname{Rest}_T$. Define $\Delta_T(s/T, u) := c_{X,T}(s/T, u) - c_X(s/T, u)$. We first show that $\operatorname{TV}(\Delta_T(\cdot, u))$ is uniformly bounded in $u$. For all $I \in \{1, \ldots, T\}$ and every sequence $0 < a_1 < a_2 < \cdots < a_I < 1$, we can write

$$\Delta_T(a_i, u) - \Delta_T(a_{i-1}, u)$$
$$= \sum_{j=-\infty}^{-1} \sum_{k=-\infty}^{\infty} \left\{ S_j\left(\frac{k}{T}\right) - S_j(a_i) \right\} \psi_{jk}([a_i T]) \psi_{jk}([a_i T] + u)$$
$$- \sum_{j=-\infty}^{-1} \sum_{k=-\infty}^{\infty} \left\{ S_j\left(\frac{k}{T}\right) - S_j(a_{i-1}) \right\}$$
$$\times \psi_{jk}([a_{i-1} T]) \psi_{jk}([a_{i-1} T] + u) + O(T^{-1}),$$

where the $O(T^{-1})$ term comes from the approximation (2.2). Now, replace $k$ by $k + [a_i T]$ in the first sum and by $k + [a_{i-1} T]$ in the second one. The main term becomes

$$\sum_{j=-\infty}^{-1} \sum_{k=-\infty}^{\infty} \left\{ S_j\left(\frac{k}{T} + a_i\right) - S_j\left(\frac{k}{T} + a_{i-1}\right) + S_j(a_{i-1}) - S_j(a_i) \right\} \psi_{jk}(0) \psi_{jk}(u).$$

Consequently, using the Cauchy–Schwarz inequality and Definition 1,

$$\sum_{i=1}^{I} \{\Delta_T(a_i, u) - \Delta_T(a_{i-1}, u)\}$$
$$\leq 2 \sum_{j=-\log_2 T}^{-1} L_j \sum_{k=-\infty}^{\infty} |\psi_{jk}(0) \psi_{jk}(u)| + O(IT^{-1})$$
$$\leq 2\rho + K,$$

where $K$ is a constant (because $I \leq T$), leading to $\mathrm{TV}(\Delta_T(\cdot, u)) \leq 2\rho + K$, uniformly in $u$. We can now bound $\mathrm{Rest}_T$ in (B.26) as follows:

$$\mathrm{Rest}_T = 2\gamma_0 |\mathcal{R}T|^{-1} \sum_{s,u=-\infty}^{\infty} \Delta_T\left(\frac{s}{T}, u\right) \sum_{\ell=-\log_2 T}^{-1} 2^{\ell/2} \Psi_\ell(u)$$

$$= \frac{2\gamma_0}{|\mathcal{R}|} \sum_{s,u} \int_{s/T}^{(s+1)/T} dz \left\{ \Delta_T(z, u) + \Delta_T\left(\frac{s}{T}, u\right) - \Delta_T(z, u) \right\} \sum_\ell 2^{\ell/2} \Psi_\ell(u)$$

$$\leq \frac{2\gamma_0}{|\mathcal{R}|} \int_0^1 dz \sum_u |\Delta_T(z, u)|$$

$$+ \frac{2\gamma_0}{|\mathcal{R}|} \sum_{s,u} \int_0^{1/T} dz \left| \Delta_T\left(\frac{s}{T}, u\right) - \Delta_T\left(z + \frac{s}{T}, u\right) \right|,$$

as $|\Psi_\ell(u)|$ is uniformly bounded by 1. From Proposition 1, the first term is $O(|\mathcal{R}T|^{-1})$. Using (B.5) and the fact that $\mathrm{TV}(\Delta_T(\cdot, u))$ is uniformly bounded in $u$, the second term is also $O(|\mathcal{R}T|^{-1})$.

In (B.26), we now expand $c_X(s/T, u)$ using (2.6). By the definition of the matrix $A$, we get

$$\mathrm{tr}(M^\star_{j,\mathcal{R};T} - M_{j,\mathcal{R};T}) \geq |\mathcal{R}T|^{-1} \sum_s \sum_m S_m\left(\frac{s}{T}\right) \sum_\ell (2\gamma_0 - 2^{-\ell/2} A^{-1}_{j\ell}) 2^{\ell/2} A_{m\ell}$$

for $T$ large enough. The lower bound is derived from the definition of $\gamma_0, \gamma_1$ and Assumption 2. The upper bound is derived using $\mathrm{tr}(M^\star_{j,\mathcal{R};T} - M_{j,\mathcal{R};T}) \leq \mathrm{tr}(M^\star_{j,\mathcal{R};T})$ from (B.25), Assumption 1 and the fact that $|\Psi_\ell(u)| \leq 1$ uniformly in $\ell < 0$ and $u \in \mathbb{Z}$.

2. We prove (B.23). The first inequality is (B.1). From (B.4), we write $\|M^\star_{j,\mathcal{R};T} - M_{j,\mathcal{R};T}\|_2^2 \leq \|\Sigma^{1/2}\|_{\mathrm{spec}}^4 \|U^\star_{j,\mathcal{R};T} - U_{j,\mathcal{R};T}\|_2^2$. Then, using Lemma B.2, (A.5) and $\sqrt{\mathcal{L}_\ell \mathcal{L}_m} \leq 2^{-(\ell+m)/2} \mathcal{L}_{-1}$,

$$\tfrac{1}{2} \|U^\star_{j,\mathcal{R};T} - U_{j,\mathcal{R};T}\|_2^2$$

$$\leq \|U^\star_{j,\mathcal{R};T}\|_2^2 + \|U_{j,\mathcal{R};T}\|_2^2$$

$$\leq 4\gamma_0^2 |\mathcal{R}|^{-2} T^{-1} \sum_{m,\ell=-\log_2 T}^{-1} 2^{(\ell+m)/2} A_{\ell m} + K_2^2 2^j |\mathcal{R}|^{-2} T^{-1}$$

$$\leq 4\mathcal{L}_{-1} \gamma_0^2 |\mathcal{R}|^{-2} T^{-1} \log_2^2(T) + O(T^{-1}).$$

The result follows from Lemma B.3.

3. We prove (B.24). For $T$ large enough, $\lambda_T$ is strictly positive. Using Proposition B.1 and defining $p_T^2 = \mathrm{Var}(\underline{Z}'_T (M^\star_{j,\mathcal{R};T} - M_{j,\mathcal{R};T}) \underline{Z}_T) = 2\|M^\star_{j,\mathcal{R};T} -$

$M_{j,\mathcal{R};T}\|_2^2$ and $q_T = \|M_{j,\mathcal{R};T}^\star - M_{j,\mathcal{R};T}\|_{\mathrm{spec}}$, we can write

$$\Pr(\underline{Z}_T'(M_{j,\mathcal{R};T}^\star - M_{j,\mathcal{R};T})\underline{Z}_T > \lambda_T)$$

$$\leq \exp\left(-\frac{1}{2} \cdot \frac{(\mathrm{tr}\, M_{j,\mathcal{R};T})^2}{p_T^2 \log_2^2 T + 2q_T \mathrm{tr}(M_{j,\mathcal{R};T}) \log_2 T}\right).$$

(B.23) gives the rates for $p_T$ and $q_T$, leading to the result. $\square$

PROOF OF PROPOSITION 6. By Proposition 2, we have $\theta_T = Q_{j,\mathcal{R}_0} - Q_{j,\mathcal{R}_1} + O(2^{j/2}/\{\sqrt{T}(|\mathcal{R}_0| \wedge |\mathcal{R}_1|)\})$. This shows that the sign of $\theta_T$ is determined by the sign of $(Q_{j,\mathcal{R}_0} - Q_{j,\mathcal{R}_1})$ for $T$ large enough. We then consider the two cases $\theta_T > 0$ and $\theta_T < 0$.

If $\theta_T > 0$, define $\mu_T = \mathrm{E}(Q_{j,\mathcal{R}_0;T} - Q_{j,\mathcal{R};T}) > 0$ and $\lambda_T = \mathrm{tr}(M_{j,\mathcal{R}_0;T}^\star - M_{j,\mathcal{R};T}^\star) - \mu_T(1 - 1/\log_2 T)$, where the matrices $M^\star$ are defined as in Lemma B.9. Define the random set $\mathcal{P}_T = \{\underline{Z}_T'(M_{j,\mathcal{R}_0;T}^\star - M_{j,\mathcal{R};T}^\star - M_{j,\mathcal{R}_0;T} + M_{j,\mathcal{R};T})\underline{Z}_T \leq \lambda_T\}$, where $\underline{Z}_T = (Z_1, \ldots, Z_T)'$ is a vector of i.i.d. Gaussian random variables. As for the derivation of (B.24), we can use Proposition B.1 to derive

$$\Pr(\mathcal{P}_T^c) = O\left(\exp\left\{-\frac{\mu_T \sqrt{T}}{\log_2^2 T}\left(\frac{1}{|\mathcal{R}_0|^2} + \frac{1}{|\mathcal{R}_1|^2}\right)^{-1/2}\right\}\right).$$

Using decomposition (B.7) and by conditioning on $\mathcal{P}_T$,

$\Pr(\mathcal{R}$ is not rejected$|\mathcal{P}_T)$

$$\leq \Pr\{\underline{Z}_T'(M_{j,\mathcal{R}_0;T}^\star - M_{j,\mathcal{R};T}^\star)\underline{Z}_T + q_{j,\mathcal{R}_0;T}^\circ - q_{j,\mathcal{R};T}^\circ$$

$$\leq 2\eta(\sigma_{j,\mathcal{R}_0,T} + \sigma_{j,\mathcal{R},T})k_T + \lambda_T |\mathcal{P}_T\}.$$

Note that the first inequality of Proposition B.1 implies that $\Pr\{|q_{j,\mathcal{R}_0;T}^\circ - q_{j,\mathcal{R};T}^\circ| > (\sigma_{j,\mathcal{R}_0;T} + \sigma_{j,\mathcal{R};T})\lambda k_T\} \leq 2\exp(-\lambda^2 k_T^2/2)$. Therefore, by the definition of $\eta$,

$\Pr(\mathcal{R}$ is not rejected$|\mathcal{P}_T)$

$$\leq O(T^{-1}) + \Pr\{\underline{Z}_T'(M_{j,\mathcal{R}_0;T}^\star - M_{j,\mathcal{R};T}^\star)\underline{Z}_T$$

$$\leq \eta(\sigma_{j,\mathcal{R}_0,T} + \sigma_{j,\mathcal{R},T})k_T + \lambda_T |\mathcal{P}_T\}.$$

Lemma B.8 can now be used to bound this probability because $M_{j,\mathcal{R}_0;T}^\star - M_{j,\mathcal{R};T}^\star$ is a positive definite matrix and $\eta(\sigma_{j,\mathcal{R}_0,T} + \sigma_{j,\mathcal{R},T})k_T + \lambda_T \leq \mathrm{tr}(M_{j,\mathcal{R}_0;T}^\star - M_{j,\mathcal{R};T}^\star)$ for $T$ large enough. This leads to the rate $O(-\frac{\mu_T^2 T}{\log_2^2 T}(\frac{1}{|\mathcal{R}_0|^2} + \frac{1}{|\mathcal{R}|^2})^{-1})$.

If $\theta_T < 0$, then we apply the same reasoning with $\mu_T = \mathrm{E}(Q_{j,\mathcal{R}_1;T} - Q_{j,\mathcal{R};T})$ and $\lambda_T = \mathrm{tr}(M_{j,\mathcal{R}_1;T}^\star - M_{j,\mathcal{R};T}^\star) + \mu_T(1 - 1/\log_2 T)$. The result follows after the addition of all terms. $\square$

**Acknowledgments.** We are grateful to I. Gijbels, G. P. Nason, V. Spokoiny and S. Subba Rao for helpful comments and to an Associate Editor and two referees for their constructive criticism which lead to significant improvements in the paper.

Université catholique de Louvain
Institut de statistique
Voie du Roman Pays 20
B-1348 Louvain-la-Neuve
Belgium
E-mail: vanbellegem@stat.ucl.ac.be
       vonsachs@stat.ucl.ac.be